\documentclass[preprint]{elsarticle}

\usepackage{lineno,hyperref}
\usepackage[centertags]{amsmath}
\usepackage{amsfonts}
\usepackage{amssymb}
\usepackage{amsmath, bm}
\usepackage[ruled]{algorithm2e}
\usepackage{enumitem}
\usepackage{multicol}
\usepackage{multirow}
\usepackage{makecell}
\usepackage{ragged2e}
\usepackage{color}
\renewcommand{\justifying}{\leftskip=0pt \rightskip=3pt plus 0cm}
\renewcommand{\justify}{\leftskip=0pt \rightskip=0pt}

\newtheorem{remark}{Remark}
\journal{}
\SetKwBlock{pre}{Pre-training step:}{end}
\SetKwBlock{formal}{Formal training step:}{end}
\numberwithin{table}{section}

\numberwithin{equation}{section}
\numberwithin{figure}{section}
\begin{document}
\begin{frontmatter}
\title{Pre-training strategy for solving evolution equations based on physics-informed neural networks}

\author[1,2]{Jiawei Guo}
\author[2]{Yanzhong Yao\corref{cor1}}
\ead{yao\_yanzhong@iapcm.ac.cn}
\author[2,3]{Han Wang}
\author[2]{Tongxiang Gu\corref{cor1}}
\ead{txgu@iapcm.ac.cn}

\cortext[cor1]{Corresponding author}
\address[1]{Graduate School of China Academy of Engineering Physics, Beijing 100088, China}
\address[2]{Laboratory of Computational Physics, Institute of Applied Physics and
Computational Mathematics, Beijing 100088, China}
\address[3]{HEDPS, CAPT, College of Engineering, Peking University, Beijing, China}

\begin{abstract}
The physics informed neural network (PINN) is a promising method for solving time-evolution partial differential equations (PDEs). 
However, the standard PINN method may fail to solve the PDEs with strongly nonlinear characteristics or those with high-frequency solutions.
To address this problem, we propose a novel method named pre-training PINN (PT-PINN) which can improve the convergence and accuracy of the standard PINN method by combining with the resampling strategy and the existing optimizer combination technique.
The PT-PINN method transforms the difficult problem on the entire time domain to relatively simple problems defined on small subdomains. 
The neural network trained on small subdomains provides the neural network initialization and extra supervised learning data for the problems on larger subdomains or on the entire time-domain. 
By numerical experiments, we demonstrate that the PT-PINN succeeds in solving the evolution PDEs with strong non-linearity and/or high frequency solutions, including the strongly nonlinear heat equation, the Allen-Cahn equation, the convection equation with high-frequency solutions and so on, and that the convergence and accuracy of the PT-PINN is superior to the standard PINN method.
The PT-PINN method is a competitive method for solving the time-evolution PDEs.
\end{abstract}


\end{frontmatter}


\section{Introduction}

Many practical applications such as radiation energy transfer problems, convection problems and phase separation problems require investigations of time-evolution partial differential equations (PDEs).
Traditional numerical methods, like the finite difference, the finite volume and the finite element schemes, face many challenges in solving these equations, and many researchers have devoted their energies to constructing high confidence discrete schemes and high efficiency iteration methods.
The traditional numerical methods face the following difficulties:
(1) Generally, the design of discrete schemes is closely related to dimensionality, and most traditional numerical algorithms face great difficulties while the computational model changes from two dimensions to three dimensions. 
(2) For the evolution equation with strong nonlinear characteristics, it is necessary to study the corresponding nonlinear iterative algorithms according to the type of the equation. Moreover, these algorithms are not of generality.

In recent years, machine learning techniques for solving PDEs have developed vigorously.
In 1998, Lagaris et al.~first tried to solve the PDEs by using machine learning methods.
They expressed the solution with artificial neural networks (ANNs). 
By designing a reasonable loss function, the output of the network met the equation and definite conditions \cite{Lagaris1}. 
In a further work, they modified the network architecture to exactly satisfy the initial conditions and Dirichlet boundary conditions for the case of complex geometric boundary, which improved the prediction performance \cite{Lagaris2}.

\begin{figure}[]
\centering
\includegraphics[width=0.6\textwidth]{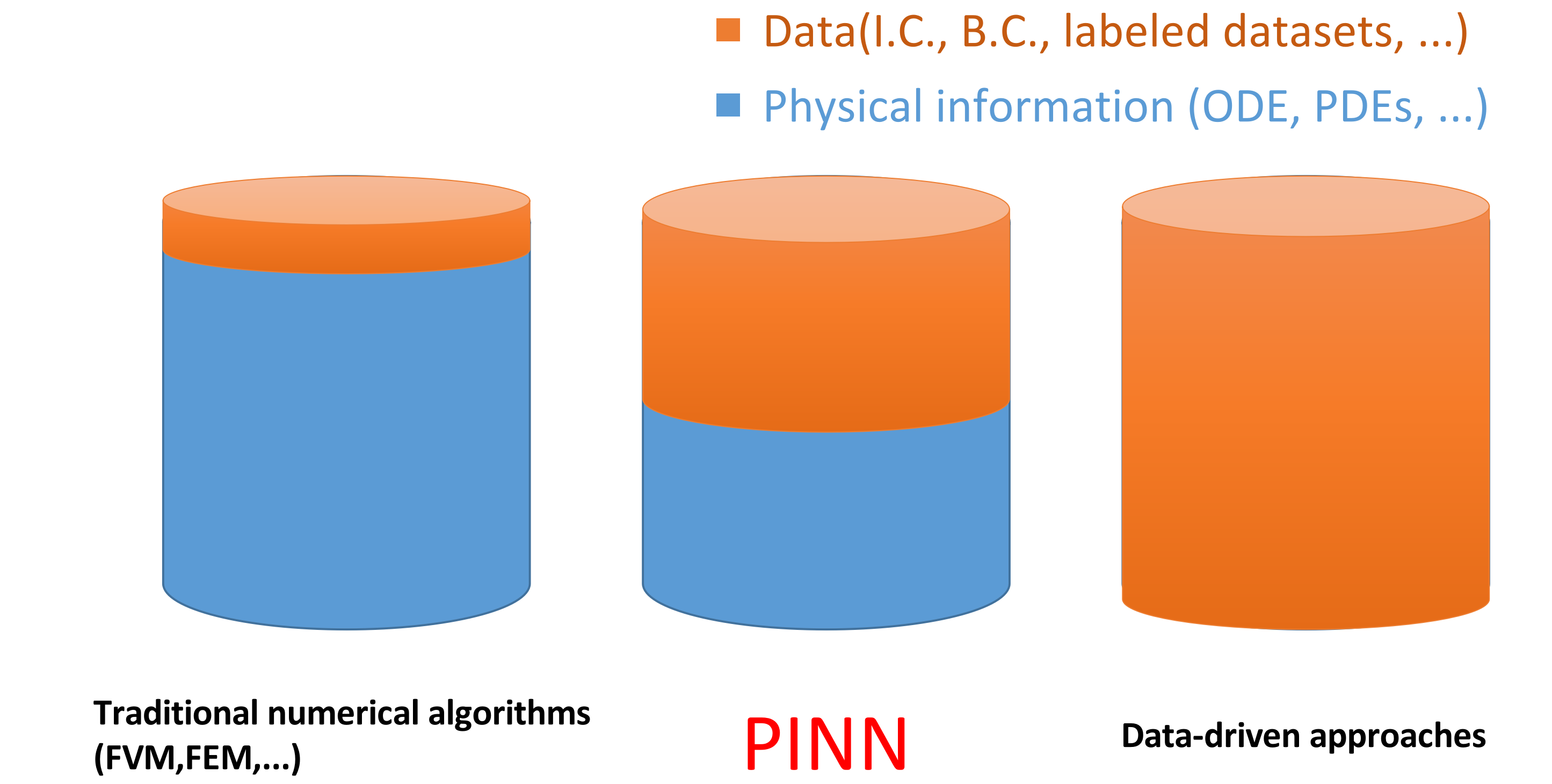}
\caption{Data characteristics of traditional numerical algorithms, PINN and data-driven approaches.}
\label{Fig1.1}
\end{figure}

With the development of automatic differentiation technique, Raissi et al.~\cite{RAISSI1,Karniadakis} revisited these methods by using modern computational tools under the named as physics informed neural networks (PINN) and applied them to solve various PDEs. 
PINN is a deep learning framework for solving forward and inverse problems involving nonlinear PDEs, which was first clearly proposed in 2019~\cite{RAISSI1}. 
By introducing the information of equations into the training process as a part of the training goal, and by using the automatic differentiation technique~\cite{Baydin}, PINN provides basically accurate and physically interpretable predictions~\cite{Karniadakis}. 
Compared with traditional numerical algorithms and common data-driven approaches, PINN has the characteristics as shown in Figure \ref{Fig1.1} in the aspects of using both data and physical information.  
In the Figure, the ``physical information'' stands for the physical law expressed by PDEs, while the ``data'' represents the supervised learning data named as labeled data, which includes the initial condition (I.C.), the boundary condition (B.C.) and other auxiliary data. 
Traditional numerical algorithms construct the discrete schemes on the computational grid based on the definite conditions and the equation, and then use the corresponding solution methods to obtain the approximate numerical solutions.
Commonly used data-driven approaches directly train the neural networks by using a large amount of labeled data (also named as the training set), so that they obtain the network prediction on the test set which is independent from the training set.
The PINN seamlessly integrates the physical laws and the labeled data through the loss function. 
The physical laws are reflected by the residuals of the PDE at sample points named as residual points. 
The labeled data are obtained from definite conditions like I.C.~and B.C.~of the equation. 
The neural network can be obtained by optimizing the loss function, and if the optimization succeeds the neural network prediction obeys the physical laws and labeled data.

Compared to the traditional numerical methods for solving PDEs, PINN has three virtues: 
(1) PINN is a mesh free method, and it completely avoids the difficulties, like the serious distortion of grids, arising in the traditional numerical algorithms using Lagrangian grids;
(2) PINN avoids the curse of dimensionality \cite{Poggio,Grohs} encountered in traditional numerical calculation to some extent. 
The implementation of PINN is basically the same for one-dimensional to three-dimensional or even high-dimensional problems; 
(3) The final training result of PINN is a continuous mapping, which gives the prediction of the equation at any time and position in the solution domain.

In recent years, researchers have applied the PINN to many fields. Jin et al. explored the effectiveness of PINN to directly simulate Navier-Stokes equation (NSFnets), and used the transfer learning (TL) method to accelerate the solution~\cite{Jin}. 
Their numerical results indicate that PINN is better than traditional numerical algorithms at solving the ill-posed or inverse problems.
Mishra et al. proposed a novel machine learning algorithm based on PINN to solve forward and inverse problems for the radiation transport model.
Furthermore, they presented extensive experiments and theoretical error estimates to demonstrate that PINN is a robust method \cite{Mishra}.
In Ref.~\cite{Zobeiry}, PINN was developed to solve a linear heat transport PDE with convective boundary conditions.
Their work shows that the trained network has near real time simulation capability of problems with given boundary conditions. 
In addition, PINN has also been successfully applied to stochastic groundwater flow analysis in heterogeneous aquifers \cite{Guo}. For nonlinear conservation laws,
Jagtap et al. proposed a conservative physics-informed neural networks (cPINN) on discrete domains, which has the flexibility of domain decomposition into several non-overlapping sub-domains~\cite{Jagtap1}. 
As a development version of cPINN, the extended PINN (XPINN) offered both space and time parallelization and reduced the training cost more effectively~\cite{Jagtap2}.

PINN can be generalized to handle various types of differential equations. 
Pang et al. extended PINN to fractional PINN (fPINN) for solving fractional differential equations (FDEs) and systematically demonstrated the convergence and efficiency of the algorithm \cite{Pang}.
PINN was used to solve the ordinary differential equations (ODEs), the integro-differential equations (IDEs) and stochastic differential equations (SDEs) \cite{Ji,Yang,Zhang,Nabian}.

PINN 
may fail in training and to produce accurate predictions, when the PDE solutions contain high-frequency or multi-scale features \cite{Raissi4,Fuks}.
Mattey et al. demonstrated that the accuracy of the PINN suffers in the presence of strong non-linearity and higher order partial differential operators \cite{Mattey}. 
The loss function of PINN generally consists of multiple terms.
During the training process, different terms of the loss function compete with each other, and it may be difficult to minimize all the terms simultaneously.
Therefore, one may need to develop more efficient network architectures and optimization algorithms for different equations. 
In Ref.~\cite{Ji}, Ji et al. investigated the performance of PINN in solving stiff ODEs, and the results showed that the stiffness could be the main reason of training failure. 
Wang et al. analyzed the training dynamics of PINN through the limiting neural tangent kernel and found a remarkable discrepancy in the convergence rate of the different loss components contributing to the total training error in Refs.~\cite{Wang1,Wang2,Wang3}. 
They mitigated this problem by designing a more appropriate loss function and new optimization algorithms.
PINN is a strong form approach, and it results in a smooth prediction which may be biased from the real solution of the problems with low smoothness. 
Consequently, some methods based on the weak formulation of PDEs have been proposed.
Zang et al. developed a approach named as weak adversarial network (WAN)~\cite{Zang}. 
Kharazmi et al. adopted the weak formulation of the PDE and hp-refinement to enhance the approximation capability of networks~\cite{Karumuri}. 

At present, the theoretical analysis of PINN method is still rare. In this regard, Shin et al.~made a theoretical analysis on the solution of PDEs by PINN~\cite{Shin}, and proposed that the solution of PINN could converge to the solution of the equation under certain conditions.
Mishara et al.~gave an abstract analysis of the generalization error source of PINN in solving PDEs~\cite{Mishra2}.

It is worth mentioning that the concept of time-space variable is blurred in some literature about the continuous-time PINN method \cite{Lagaris1,Jagtap2,Lu}. 
As for the evolution PDEs, a good prediction near the initial time is the basis for the effectiveness of the PINN, and the lack of time-dependent characteristics may lead to training failure \cite{Mattey,Zang,Wight}. 
Recently, some time-machine methods are proposed to solve a type of evolution equations \cite{Mattey,Wight,Krishnapriyan,Haitsiukevich}.
If the evolution equations have strong non-linearity or high-frequency solution, the training of standard PINN method struggles, and this is the topic that we are focusing on.

In the present paper, we investigate the solution method of the evolution PDEs by PINN and analyse why PINN fails to train in some cases. 
In order to improve the accuracy of prediction and the robustness of training algorithms, a novel PINN method is proposed based on the \textit{pre-training (PT) strategy}.
We name the new method as \textit{pre-training PINN (PT-PINN)}. 
Moreover, we adopt the residual points resampling technique \cite{Lu} and the optimizers combination strategy \cite{He}, which are very helpful for improving the accuracy and efficiency of the training.  

The rest of this paper is organized as follows: In section 2, a typical evolution PDE and its standard PINN method are introduced. 
In section 3, the pre-training, the optimizers combination and  the resampling strategies for improving the standard PINN are presented. The difference between the PT-PINN and the related work is also discussed.
In section 4, some examples are provided to demonstrate the performances of the PT-PINN method.
In the last section, we draw a conclusion about our work.

\section{Governing equation and the standard PINN method}

In this section, we first introduce a typical evolution PDE, and then  its standard PINN method. 
Moreover, we discuss the issues of the standard PINN.

\subsection{Physical model and the governing equation}

Evolution PDEs refer to the general name of many important PDEs with the temporal variable $t$, mainly including heat equations, convection equations, reaction-diffusion equations, the Schr\"{o}dinger equation, hydrodynamic equations, the Korteweg-de Vries (KdV) equation and so on. 

In this paper we focus on the evolution PDEs of the form
\begin{equation}\label{21}
u_t + \mathcal{N}[u]=0,\quad X\in\Omega,\,t\in(0,T],
\end{equation}
with the I.C.~and the Dirichlet B.C.,
\begin{equation}\label{eq:ic-and-bc}
\begin{aligned}
&u(X,0)=\mathcal{I}(X),\quad X\in\Omega,\\
&u(X,t)=\mathcal{B}(X,t),\quad X\in\partial\Omega,\,t\in(0,T].
\end{aligned}
\end{equation}
In Eqs.~\eqref{21} and \eqref{eq:ic-and-bc}, $u(X,t)$ denotes the solution with spatial variable $X$ and domain $\Omega$, and $\mathcal{N}$ represents a partial differential operator with respect to the spatial variable.
$\mathcal{I}(X)$ and $\mathcal{B}(X,t)$ are bounded in their domain. 
We suppose the given B.C. is well-posed. For different models, the B.C. can also be other types, such as 
 Neumann, Robin or periodic types.
We focus our attention on the PDE that has a continuous solution on the whole domain.

An example of the evolution equation is the heat equation defined by
\begin{align}
   \mathcal N[u] = - \nabla\cdot(\kappa(u,X,t)\nabla u) - Q(u, X, t)
\end{align}
where $u$ is the temperature function to be solved,  $\kappa$ is the continuous heat conduction coefficient dependent on $u$, and $Q$ is the source term.

\subsection{Standard PINN}
\label{sec:std-pinn}

\begin{figure}[]
\centering
\includegraphics[width=0.7\textwidth]{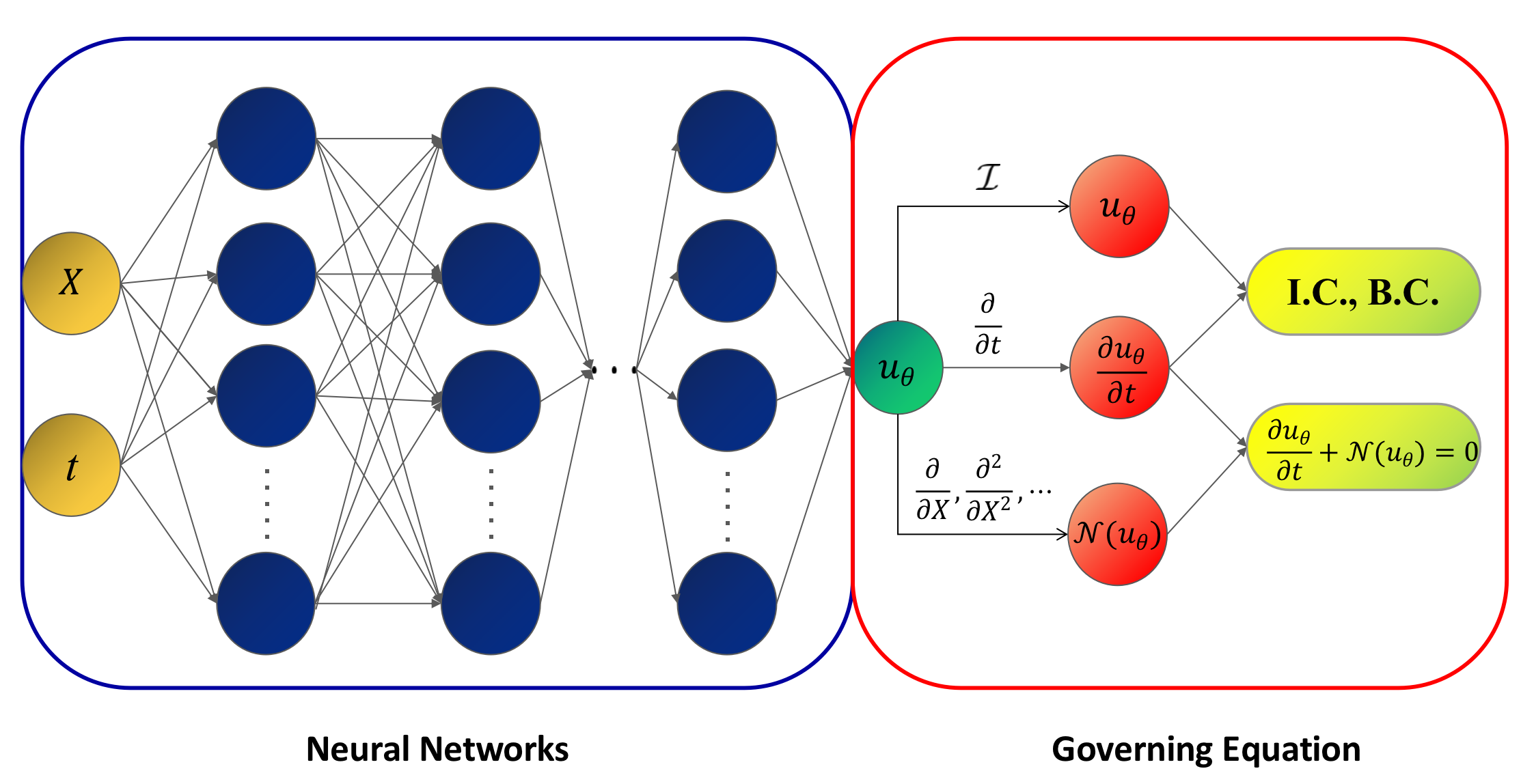}
\caption{Schematic plot of the architecture of PINN for Eq.~\eqref{21}}
\label{Fig2.1}
\end{figure}

By the PINN method, the solution of Eq.~\eqref{21} is approximated by a deep neural network (DNN), which takes spatial ($X$) and temporal ($t$) variables as inputs, and outputs the approximated solution $u_\theta$, where $\theta$ denotes the {\itshape neural network parameters} including the weights and biases of neural networks. 
Figure~\ref{Fig2.1} schematically shows the architecture of the PINN. 
The partial differential operator $\mathcal N$ and $\frac{\partial}{\partial t}$ are implemented by using automatic differentiation, which can be easily realized in the deep learning framework like the Tensorflow \cite{TF} or PyTorch \cite{PyTorch}.
The neural network parameters $\theta$ are obtained by optimizing the loss function which is composed of three parts as follows:
\begin{align}
&\mathcal{L}(\theta;\Sigma)=w_i\mathcal{L}_i(\theta;\tau_i)+w_b\mathcal{L}_b(\theta;\tau_b)+w_r\mathcal{L}_r(\theta;\tau_r),\label{loss}\\
&\mathcal{L}_i(\theta;\tau_i)=\frac{1}{N_i}\sum_{i=1}^{N_i}\left|u_\theta(X_i,0)- \mathcal{I}(X_i)\right|^2,\label{loss1} \\
&\mathcal{L}_b(\theta;\tau_b) = \frac{1}{N_b}\sum_{i=1}^{N_b}\left|u_\theta\left(X_i,t_i\right)-\mathcal{B}(X_i,t_i)\right|^2,\label{loss2}\\
&\mathcal{L}_r(\theta;\tau_r)=\frac{1}{N_r}\sum_{i=1}^{N_r}\left|\frac{\partial u_\theta}{\partial t}\left(X_i,t_i\right)
+ \mathcal N[u]\left(X_i,t_i\right)
\right|^2.\label{loss3}
\end{align}
In Eqs. (\ref{loss}-\ref{loss3}), 
$w_i$, $w_b$ and $w_r$ are the weights for three parts of the loss function.
$\mathcal{L}_i$ and $\mathcal{L}_b$ respectively denote the I.C.~and B.C.~supervised learning losses.
The labeled data sets for training the I.C. and B.C. are denoted by
$$\begin{aligned}
&\tau_i=\{\left(X_i,\mathcal{I}(X_i)\right)|X_i\in\Omega\}_{i=1}^{N_i},\\
&\tau_b=\{\left(X_i,t_i,\mathcal{B}(X_i,t_i)\right)|(X_i,t_i)\in\partial\Omega\times(0,T]\}_{i=1}^{N_b}.
\end{aligned}
$$
By training the neural network in the labeled data sets, the prediction of PINN is close to the real solution $u$ in sets $\tau_i$ and $\tau_b$.
In Eqs.~\eqref{loss} and \eqref{loss3}, $\mathcal{L}_r$ represents the residual loss defined in the residual data set 
$$\tau_r=\{\left(X_i,t_i\right)\in\Omega\times(0,T]\}_{i=1}^{N_r}.$$
By minimizing Eq.~\eqref{loss3} in $\tau_r$, PINN provides the prediction approximately satisfying the governing equation Eq.~\eqref{21} in domain $\Omega \times (0, T]$. 
$N_i, N_b$ and $N_r$ denotes the size of I.C.~($\tau_i$), B.C.~($\tau_b$) and the residual ($\tau_r$) data sets, respectively. 
$\Sigma$ denotes the collection of them, i.e.~$\Sigma = \{\tau_i,\tau_b,\tau_r\}$.
The training of PINN is a procedure of solving the optimization problem
\begin{equation}\label{MIN}
\bar{\theta} = \arg\min_\theta\mathcal{L}(\theta;\Sigma).
\end{equation}

We use the training error $\left\|u_\theta-u\right\|$ to measure the performance of the PINN method. The $L_2$-norm relative error is defined as follows
\begin{equation}\label{r2}
\left\|\varepsilon\right\|_2=\frac{\sqrt{\sum_{i=1}^N\left|u_\theta(X_i,t_i)-u(X_i,t_i)\right|^2}}{\sqrt{\sum_{i=1}^N\left|u(X_i,t_i)\right|^2}},
\end{equation}
where $u(X_i,t_i)$ is the real solution or the reference solution, and $u_\theta(X_i,t_i)$ is the neural network prediction for a point lying in a test set $\{(X_i, t_i)\in\Omega\times(0,T])\}^{N}_{i=1}$. We construct the test set by randomly sampling points from the spatial-temporal domain, and the size of test set is $N=10^4$ for all numerical examples presented in Sect.~\ref{sec:examples}.
In addition, the accuracy of the trained neural network is assessed by taking the $L_1$-norm and $L_\infty$-norm absolute errors, which are defined as follows:
\begin{equation}
\left\|e\right\|_1=\frac{1}{N}\sum_{i=1}^N\left|u_\theta(X_i,t_i)-u(X_i,t_i)\right|,
\end{equation}
\begin{equation}
\left\|e\right\|_\infty=\max_{1\leq i\leq N}\left|u_\theta(X_i,t_i)-u(X_i,t_i)\right|.
\end{equation}

\subsection{Several issues of PINN for solving evolution PDEs}
\label{sec:issues}

Despite the great achievement of the PINN method, it still has the following issues:
\begin{itemize}
\item The PINN loss includes the supervised learning parts and the residual learning part as the learning objectives. Since the dimensions of the parts are different, general stochastic gradient optimization algorithms may fail to train PDEs with solutions of large gradients and/or high-frequencies~\cite{Wang1,Moseley}.
\item For the strongly nonlinear evolution PDEs, the loss function is more complex, and the training of the standard PINN is sometimes trapped by the local minima, which leads to the failure of training~\cite{Mattey,Wight,Wang4}.
\item The PINN introduces the information of the PDE by the residual learning part. 
The number and the distribution of residual data points directly affect the quality of the training. 
PINN usually requires a large number of residual data set.
However, large numbers of residual points lead to a significant increase in the computational cost~\cite{Lu}.
\item 
The prediction provided by PINN method tends to enforce smoothness and to violate the I.C.~\cite{Zang}.
A poor initial guess of the optimizer easily results in a homogeneous solution, see the results of the standard PINN method in Sect.~\ref{reaction} and Sect.~\ref{S44}.
This issue is particularly prominent in the problems with periodic boundary conditions~\cite{Krishnapriyan,Haitsiukevich}.
\end{itemize}

These problems motivate the main goal of this paper, and some strategies are proposed to mitigate them in the subsequent section.

\section{The PT-PINN method for evolution PDEs}

In this section we introduce the idea of pre-training and the optimizer combination and the resampling techniques to improve the  reliability and the convergence of the standard PINN method.

\subsection{The pre-training strategy}\label{3.1}

The pre-training strategy is composed of two steps:
\begin{enumerate}
\item The pre-training step. 
Eq.~\eqref{21} is trained in the time interval $[0,T_p]$ by using the standard PINN method as the pre-training step. 
We call $[0,T_p]$ the \textit{pre-training interval}. Note that the length of the pre-training interval is a tunable parameter and its choice is problem dependent. 
The data set $\Sigma_p = \{\tau_{ip},\tau_{bp},\tau_{rp}\}$ of the pre-training step is defined as:
$$
\begin{aligned}
&\tau_{ip}=\{\left(X_i,\mathcal{I}(X_i)\right)|X_i\in\Omega\}_{i=1}^{N_{ip}},\\
&\tau_{bp}=\{\left(X_i,t_i,\mathcal{B}(X_i,t_i)\right)|(X_i,t_i)\in\partial\Omega\times(0,T_p]\}_{i=1}^{N_{bp}},\\
&\tau_{rp}=\{\left(X_i,t_i\right)|(X_i,t_i)\in\Omega\times(0,T_p]\}_{i=1}^{N_{rp}},
\end{aligned}
$$ 
where $N_{ip},N_{bp},N_{rp}$ are the number of initial data, boundary data and residual data in the pre-training step, respectively.
The corresponding optimization problem of the pre-training step is
\begin{equation}\label{pretrain}
\theta_p=\arg\min_{\theta}\mathcal{L}(\theta;\Sigma_p).
\end{equation}

The value of $\theta$ is usually randomly initialized when solving the problem~\eqref{pretrain}.
\item The formal training step.
The additional supervising information of the solution can effectively improve the performance of PINN~\cite{He}.
The pre-trained solution from step 1 provides {extra supervised learning points} for the formal training. 
The loss function is formulated as the following form:
\begin{equation}\label{fortrain}
\mathcal{L}(\theta; \Sigma, \theta_p, \tau_{sp})= \mathcal L(\theta;\Sigma) + w_{sp}\mathcal{L}_{sp}(\theta; \theta_p, \tau_{sp}),
\end{equation}
where $w_{sp}$ is the weight of the extra supervised learning part.
$\mathcal{L}_{sp}$ is loss of the supervised learning formulated by 
\begin{align}\label{eq:loss-sp}
\mathcal{L}_{sp}(\theta; \theta_p, \tau_{sp}) = \frac{1}{N_{sp}}\sum_{i=1}^{N_{sp}}\left|u_{\theta}(x_i,t_i)- u_{\theta_p}(x_i,t_i)\right|^2,
\end{align}
with the data set $\tau_{sp}$ defined by
\begin{align}\label{eq:data-sp}
    \tau_{sp} = \{
    (X_i, t_i, u_{\theta_p}(X_i,t_i)) \, \vert \, 
    (X_i,t_i) \in \Omega\times(0,T_p]
    \}_{i=1}^{N_{sp}},
\end{align}
and $N_{sp}$ being the number of points in the data set.

Utilizing the $\theta_p$ as the initial guess of the neural network parameters,  the optimization problem
\begin{equation}\label{formaltrain}
\bar\theta=\arg\min_{\theta}\mathcal{L}(\theta;\Sigma, \theta_p, \tau_{sp}),
\end{equation}
is solved to obtain an approximate solution $u_{\bar\theta}(X,t)$ of the evolution PDE Eq.~\eqref{21} in the domain $\Omega\times(0, T]$.
\end{enumerate}

Regarding the supervised learning in the formal training step, we have the following remarks.

\begin{remark}
It is fully validate the use an empty supervised learning data set (i.e.~$N_{sp} = 0$) in the formal training step. In this case, the loss function of the formal training step reduces to the standard PINN.
\end{remark}

\begin{remark}
The idea of adding extra supervised points is closely related to the label propagation technique~\cite{Haitsiukevich}, and the extra supervised learning part is helpful to improve the training performance. For complex problems, this technique is more prominent, and the numerical examples of this work verify this phenomenon.
\end{remark}

The reasons and advantages of using the pre-training strategy is summarized as follows:
\begin{itemize}
\item 
Successful training near the initial condition is the beginning of getting meaningful results from the PINN in the entire time-domain~\cite{Wang4}.
\item 
Compared with training in the whole time domain, training in a short time segment is easier than larger time intervals \cite{Mattey,Wight}.
\item 
If the solution of the PDE is smooth, the solution obtained by the pre-training has certain extrapolation properties \cite{Karniadakis} in the time direction, which is the reason why it can be used as a good initial guess for the formal training step.
\item The pre-training step can also provide an extra supervised learning part, which is important to improve the performance of the formal training.
\end{itemize}

\subsection{Multiple pre-training strategy}

The length of pre-training interval should not be too large to ensure the accuracy of pre-training step.
However, for some difficult cases, a single pre-training step with a short pre-training interval may not provide a good enough initial guess for the formal training. 
In order to overcome this problem, we propose a \textit{multiple pre-training strategy}. 

With the multiple pre-training strategy, several pre-training  intervals are designed as follows:
\begin{align}
    (0,T_1],(0,T_2],\cdots,(0,T_{k+1}], \quad 
    \textrm{with}\ T_1 < T_2 < \dots < T_{k+1} = T.
\end{align}
In the first pre-training interval $(0,T_1]$, the standard pre-training step is reused.
The resultant neural network parameters are denoted by $\theta_p^{(1)}$.
We denote the pre-training data set of the $i$-th interval by $\Sigma_p^{(i)}$.
In any pre-training interval $(0,T_i]$, $i>1$, the following optimization problem 
\begin{align}\label{eq:mtp-formal}
    \theta_p^{(i)} = \arg\min_{\theta} \mathcal L(\theta; \Sigma_p^{(i)}, \theta_p^{(i-1)}, \tau_p^{(i-1)}),
\end{align}
is solved with the neural network parameters $\theta$ initialized with $\theta_p^{(i-1)}$.
In Eq.~\eqref{eq:mtp-formal}, the supervised learning data set $\tau_p^{(i-1)}$ is generated by using the solution pre-trained in the interval $(0, T_{i-1}]$.
The training in the last interval $(0, T_{k+1}] = (0, T]$ gives an approximation to the solution of the evolution PDE.




\begin{remark}
As for the heat equation, the diffusion operator has the smoothing characteristics.
Numerical experiments of this work show that, for the heat equation with smooth boundary conditions and weak source term, one pre-training interval ($k=1$) is enough; for the heat equation with strongly nonlinear or non-smooth source item and boundary conditions, two or more pre-training intervals ($k>1$) may be required.
\end{remark}

\subsection{The choice of optimizers}
\label{3.2}
The PINN obtains prediction of the equation by solving non-convex optimization problem. Thus, training optimizers and training strategies are crucial factors for achieving good accuracy. In this study, we chose two training optimizers: Adam \cite{Kingma} and L-BFGS \cite{Byrd} which are regularly used in training the PINN. 

Adam is a stochastic gradient descent method. 
The method is efficient when the optimization problem involves a lot of data or parameters, and it is often the first choice of the machine learning tasks.
Based on our experience,
the convergence speed of L-BFGS is faster than Adam, because L-BFGS uses the approximations to the second-order derivatives of the loss function, while Adam only relies on the first-order derivatives. 
However, as a quasi-second order optimization algorithm, L-BFGS is sensitive to the initialization, and it is easy to be trapped in bad local minima.
We thus use the Adam optimizer for a number of iterations to provide a good initial neural network, then switch to the L-BFGS optimizer until convergence. 
The loss function involved in this work, such as Eqs.~\eqref{loss}\eqref{fortrain}, are all trained by the combination of Adam and L-BFGS optimizers.



\subsection{Residual points resampling strategy}


The standard PINN generates all residual points in the data preparation stage before training, and uses all these residual points throughout the training procedure.
We instead dynamically and randomly sample the residual points from the spatial-temporal domain during the training process. 
Our numerical experiments show that resampling all residual data points in each training step could lead to a decrement of the stability and accuracy of the training. 
We thus resample a certain ratio, denoted by $\eta$ ($0\leq \eta\leq 1$) of the residual data points every $K$ training steps. 
The last resampling step that is usually smaller than the total number of Adam training steps is denoted by $F$. 
The resampling in the L-BFGS training steps usually lead to large and harmful fluctuations in the loss function.
For this reason, the resampling strategy is only adopted with the Adam optimizer.
Taking the optimization problem Eq.~\eqref{MIN} as an example, the resampling process and the optimizer choice are summarized in Algorithm~\ref{AL2.1}. 
The optimizer combination and residual data resampling  strategies of the optimization problems Eqs.~\eqref{pretrain} and \eqref{formaltrain} can be written analogously. 

If $F=N_{iter}$ and $\eta=1$, Algorithm \ref{AL2.1} becomes the resampling strategy reported in Ref.~\cite{Lu}.
For $K=F$ and $\eta=0$, Algorithm \ref{AL2.1} reduces to the standard PINN.





\begin{algorithm}[]\label{AL2.1}
        \caption{Optimizer combination and resampling strategy}
        \KwIn{\\
        {\justifying Neural network structure; Training set $\Sigma = \{\tau_i,\tau_b,\tau_r\}$;\\ Iteration number $N_{iter}$ of Adam algorithm; Resampling ratio $\eta$; \\
        Resampling interval $K$; Resampling termination step $F$. }}
        \For{$i=1,\cdots,N_{iter}$}{{\justifying
            Minimize the loss function $\mathcal{L}(\theta,\Sigma)$ by using Adam optimizer.}\\
            \If{$i\%K=0$ and $i<F$}
            {
            {
            \justifying Resample $\eta N_r$ residual points from $\Omega\times (0,T]$ and randomly replace $\eta N_r$ residual points in $\tau_r$.
            The new set is denoted by $\tau_r'$;
            }
            \\
            {
            \justifying Redefine the training data set 
            $\Sigma = \{\tau_i,\tau_b,\tau_{r'}\}.$
            }
            }
    }
    	{\justify  Minimize the loss function $\mathcal{L}(\theta;\Sigma)$ on the base of the training set $\Sigma$ until convergence by using L-BFGS optimizer.}
\end{algorithm}


\subsection{The procedure of PT-PINN method}

\begin{figure}[]
\centering
\includegraphics[width=1\textwidth]{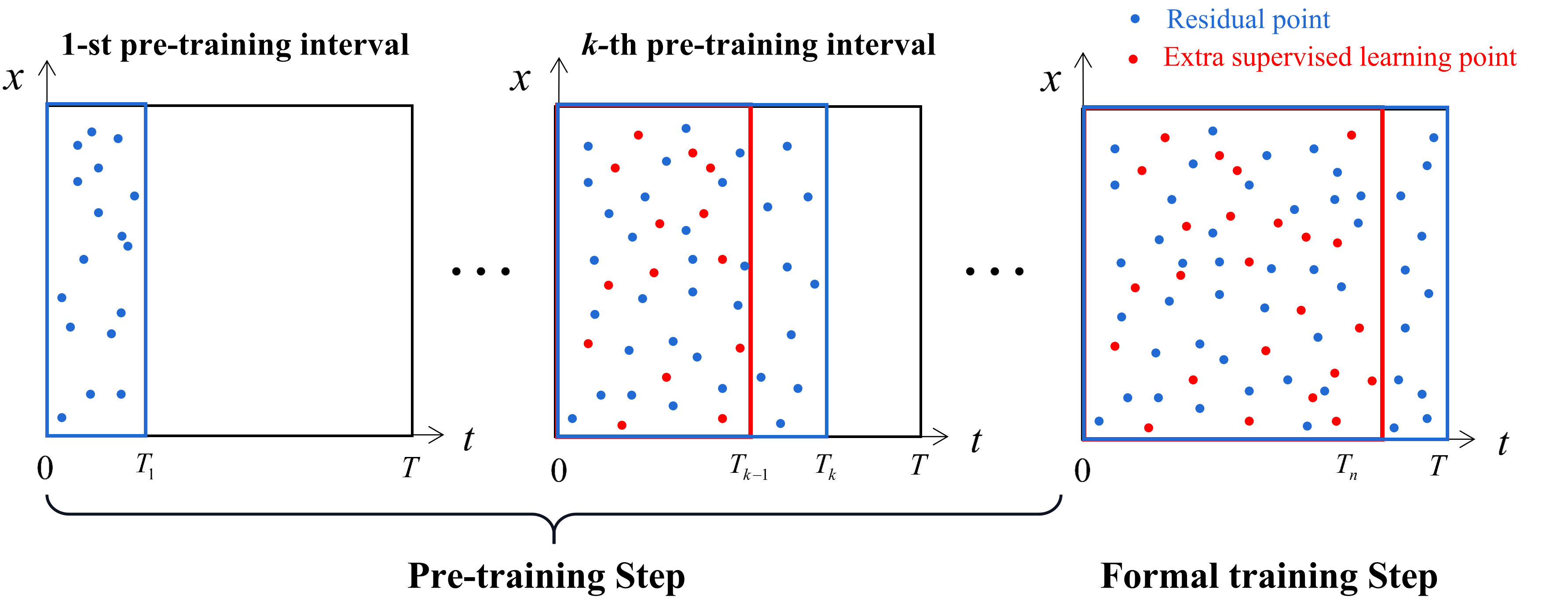}
\caption{Illustration of the PT-PINN method with $n$ pre-training intervals.
}
\label{Fig3}
\end{figure}

In this subsection, we introduce the algorithm (Algorithm~\ref{AL2.2}) of the PT-PINN for solving evolution PDEs~\eqref{21}, and the training data (residual points and extra supervised points) in each training step is shown in Figure~\ref{Fig3}. It should be noted that many details of the tunable parameters will be explained by the numerical examples in Sect.~\ref{sec:examples}. 

\begin{algorithm}[]\label{AL2.2}
        \caption{PT-PINN method for the evolution equation \eqref{21}}
        \KwIn{\\
        Neural network structure; \\
        Pre-training intervals: $(0,T_1],(0,T_2],\cdots,(0,T_k]$.
        }
        \pre{
        Neural network parameters initialization;\\ 
        {\justify Solve the standard pre-training optimization problem Eq.~\eqref{3.1} in the first interval $(0,T_1]$. Output the neural network parameters $\theta_p^{(1)}$;}\\
        \For{$i=2,\cdots,k$}{
        {\justify Generate the supervised learning data set $\tau_{sp}^{(i-1)}$ in $\Omega\times(0, T_{i-1}]$ and the training data set $\Sigma_p^{(i)}$ in $\Omega\times(0,T_i]$;}\\
        Minimize the loss function by using Algorithm \ref{AL2.1}: 
            $$\theta_p^{(i)} = \arg\min_\theta\mathcal{L}\left(\theta;\Sigma_p^{(i)},\theta_p^{(i-1)},\tau_{sp}^{(i-1)}\right);$$
        {\justify Refresh the neural network parameter $\theta_p^{(i)}$ as the initial guess of the next step.}
         }}
    	\formal{
    	{\justify 
    	Generate the supervised learning data set $\tau_{sp}^{(k)}$ in $\Omega\times(0, T_{k}]$ and the training set $\Sigma$ in $\Omega\times(0,T]$;\\
    	Use the neural network parameter $\theta_p^{(k)}$ from the last pre-training step as the initial guess and optimize Eq.~\eqref{fortrain} by using Algorithm \ref{AL2.1}}
    	$$\bar\theta=\arg\min_{\theta}\mathcal{L}(\theta;\Sigma, \theta_p^{(k)}, \tau_{sp}^{(k)}).$$
    	}
\end{algorithm}

\subsection{Related Work}

As mentioned above, the prediction of the neural network has a certain extrapolation property by assuming the smoothness of the solution and by introducing the physical information.
Therefore, the pre-training not only increases the convergence speed, but also improves the reliability of the following training.
In Ref.~\cite{Mattey}, the initial condition guided learning (ICGL) method is proposed, and the key idea is to train only the initial condition in the starting stage.
This method can be regarded as a special PT-PINN by setting  the length of the pre-training interval to 0. 
We believe that this method does not introduce physical information, and the training results only meet the initial conditions and have nothing to do with the equation.
Hence, it may not have extrapolation attribute and is less helpful for the subsequent formal training. 

It is worth pointing out that our method seems to be similar to the curriculum regularization ~\cite{Krishnapriyan} method and the time marching~\cite{Mattey,Wight,Wang4} method, however, they are different in essence. 
In Ref~\cite{Krishnapriyan}, the author devised curriculum regularization method is to warm start the training by finding a good initialization, where the loss function starts from a simple PDE regularization, and becomes progressively more complex as the neural networks get trained.
Time marching method gradually advances the training domain along the time direction, which respects the natural property of the time development equation.
The approach proposed by Ref.~\cite{Wight} trains the neural network solution of time-dependent PDEs on interval $[0, t]$, and increases the upper boundary $t$ until $T$ according to the magnitude of the residual. 
The residual points remains fixed as the time marches in Ref.~\cite{Wight},  
while the residual points are resampled in each time-interval in PT-PINN (see  Algorithm \ref{AL2.1}).
The PT-PINN is also unique in the sense that the network parameters are initialized by the parameters trained in the previous pre-training interval.
In Ref.~\cite{Mattey}, the time domain $[0,T]$ is discretized into $n$ segments as $[0=T_0,T_1],[T_1,T_2],\cdots,[T_{n-1},T_n=T]$, and the solution on the time segments are parameterized by a single neural network. 
When training the time segment $(T_{k},T_{k+1}]$, the solution predicted by the previous training at $T_{k+1}$ and in the time interval $(0,T_{k-1})$ are used to supervise the training. The network parameters of each training are initialized by an ICGL training.


Each sub-problem in the PT-PINN method can be regarded as an independent problem, and the pre-training step mainly provides a better initial guess of the neural network parameters.
The extra supervised learning part is used to improve the training performance, which can achieve a high accuracy using fewer residual points, but it is not necessary.
When the supervising data is not used, the PT-PINN avoids the accumulated error which is a common issue existing in the time-marching methods.

\section{Numerical examples}
\label{sec:examples}
In this section, we present several numerical examples including spatially 1D, 2D and 3D evolution PDEs to show the effectiveness of the PT-PINN method.
We firstly test the effect of the resampling strategy (Algorithm \ref{AL2.1}) by a 3D nonlinear heat equation in Sect. \ref{S41}. 
Other examples are used to investigate the performance of the PT-PINN method (Algorithm \ref{AL2.2}). 
In Sect.~\ref{reaction}, we successfully solve a reaction system with large reaction coefficient by using PT-PINN with a single pre-training interval, which shows that a good initial guess can significantly improve the performance of the standard PINN.
In Sect.~\ref{S42}, a heat conduction problem with high-frequency solution is well solved by the PT-PINN method, which is regarded as a challenge by the standard PINN method \cite{Krishnapriyan,Haitsiukevich,Ben}.
In Sect.~\ref{S43} and Sect.~\ref{S44}, we demonstrate the virtues of the PT-PINN for the conduction equations with strongly nonlinear diffusion coefficient and strongly nonlinear source term, respectively. 
In Sect.~\ref{Con}, a convection problem with large convection coefficient is considered, and the results show the importance of the initial guess and the extra supervised learning part $\mathcal{L}_{sp}$ (defined by Eq.~\eqref{eq:loss-sp}) in the pre-training strategy.
 
The proposed PT-PINN scheme is implemented under the deep learning framework TensorFlow version 1.5.
All the variables are of \texttt{float32} data type. 
In all examples, we use \texttt{tanh} as the activation function. 
The learning rate of the Adam optimizer exponentially decays every 50 steps at a rate of 0.98 during the training processes.
All of the parameters and stopping criteria of L-BFGS optimizer are considered as suggested in Ref.~\cite{2006Numerical}.
Before the training process, the neural network parameters are randomly initialized by using the Xavier scheme~\cite{Xavier}.


\subsection{Resampling strategy for a 3D heat equation}\label{S41}
Consider a 3D nonlinear heat equation as follows:

\begin{equation}
\begin{cases}
u_t - \nabla\cdot(u\nabla u)=Q(x,y,z,t),&(x,y,t)\in(0,1)^3\times(0,1],\\
u(x,y,z,0)=2+\sin(\pi xyz),&(x,y,z)\in(0,1)^3,\\
u(0,y,z,t)=2+\sin(5\pi t),&t\in(0,1],\\
u(1,y,z,t)=2+\sin(5\pi t+\pi yz),&t\in(0,1],\\
u(x,0,z,t)=2+\sin(5\pi t),&t\in(0,1],\\
u(x,1,z,t)=2+\sin(5\pi t+\pi xz),&t\in(0,1],\\
u(x,y,0,t)=2+\sin(5\pi t),&t\in(0,1],\\
u(x,y,1,t)=2+\sin(5\pi t+\pi xy),&t\in(0,1],
\end{cases}
\end{equation}
where the source term is
\begin{align*}
Q(x,y,z,t)=&5\pi\cos (5\pi t+\pi x y z))-  \pi^2\left(x^2y^2+x^2z^2+y^2 z^2\right)\cdot\\
   &\left(\cos (10 \pi t +2\pi x y z)-2 \sin (5\pi t+\pi x y z)\right).
\end{align*}
The exact solution is $u(x,y,z,t)=2+\sin(5\pi t+\pi xyz))$. 

For this example, we only test the impact of the resampling strategy on the standard PINN. 
We set the number of Adam iteration steps to $N_{iter}=5000$, the resampling interval to $K=200$, and the resampling termination step to $F=4000$. 
The solution is approximated by a DNN that has 6 hidden layers each with 50 neurons.
We set the numbers of supervised learning points to $N_i=100,~N_b=300$, and the number of residual points to $N_r=50,100,500,1000,2000$, respectively.

\begin{figure}[]
\centering
\includegraphics[width=0.7\textwidth]{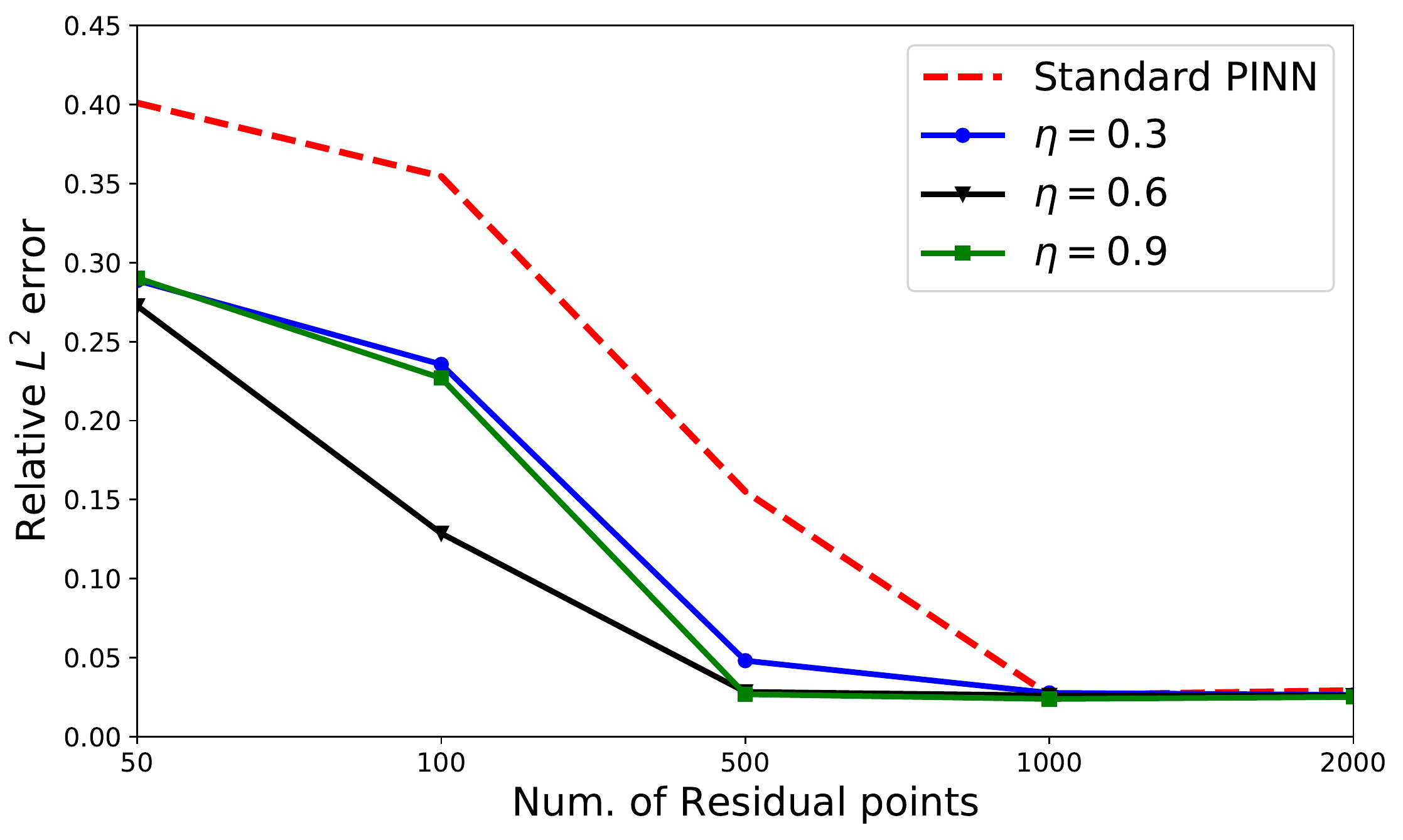}
\caption{Relative $L^2$ errors of the trainings using or not using the resampling  strategy. The resampling ratios $\eta=0.3,0.6,0.9$ are presented.
}
\label{Fig4.1}
\end{figure}

Figure \ref{Fig4.1} shows the $L^2$ relative errors for different resampling ratios $\eta$, and each result is the average of the 5 independently repeated experiments.
For the cases of fewer residual points, the resampling strategy presents obvious advantages over the standard PINN.
For the cases with relatively larger number of residual points ($N_r\geq 1000$), the resampling strategy can hardly improve the accuracy.
It can be seen that the resampling strategy using 500 residual points has the same accuracy as the standard PINN using 2000 residual points.
By using the resampling strategy, we can obtain higher accuracy with fewer residual points. 


As shown in the Figure~\ref{Fig4.1}, the optimal choice of the resampling ratio is $\eta=0.6$. 
Either a larger ($\eta=0.9$) or a smaller ($\eta=0.3$) choice will lead to an increment in the relative $L^2$ error. 
Therefore, if not stated otherwise, we choose $\eta=0.6$ in the rest numerical examples.

In Figure~\ref{Fig4.1}, only the Adam optimizer is used. 
The accuracy can be further improved by switching to the  L-BFGS optimizer after $N_{iter}$ Adam steps, as was discussed in Sect.~\ref{3.2}. 
In the case of $\eta=0.6$, the final converged relative $L^2$ error achieves $6.385\times10^{-3}$.


\subsection{Reaction equation with large reaction coefficient}\label{reaction}

The standard PINN method may fail even in solving equations of relatively simple forms. 
The reaction problem is a typical example \cite{Krishnapriyan,Haitsiukevich}. 
We consider the reaction term in the Fisher's equation as follows \cite{Krishnapriyan}:
\begin{equation}\label{rea}
\begin{cases}
\displaystyle \frac{\partial u}{\partial t} -\rho u(1-u)=0,~&(x,t)\in(0,2\pi)\times(0,1], \\
u(x,0)=h(x),~&x\in(0,2\pi),
\end{cases}
\end{equation}
where $\rho$ is the reaction coefficient. The exact solution of Eq.~\eqref{rea} is
\begin{equation}\label{eq:fisher-exact}
u(x,t)=\frac{h(x)\exp(\rho t)}{h(x)\exp(\rho t)+1-h(x)}.
\end{equation}

We compare the standard PINN method and the PT-PINN method in solving this reaction system. The I.C. is 
\begin{equation}\label{ica}
h(x)=\exp\left(-\frac{(x-\pi)^2}{\pi^2}\right).
\end{equation}
Due to the property of Eq.~\eqref{rea} and the I.C.~Eq.~\eqref{ica}, the solution should satisfy:  
\begin{equation}
u(0,t)=u(2\pi,t).
\end{equation}
We add it into the loss function as a prior knowledge, and the loss $\mathcal{L}_{b}$ is defined as follows:
\begin{align}
    \mathcal{L}_{b}(\theta;\tau_{b}) =\frac{1}{N_b}\sum_{i=1}^{N_b}\left|u_\theta(2\pi,t_i)-u_\theta(0,t_i)\right|,
\end{align}
where $\tau_b=\{(t_i\in(0,T])\}_{i=1}^{N_b}$.

It was point out that the standard PINN method can only learn this problem for small reaction coefficients~\cite{Krishnapriyan,Haitsiukevich}. 
For the case of $\rho>5$, it was shown that the standard PINN predicts a mostly homogeneous solution, which is far from the exact solution Eq.~\eqref{eq:fisher-exact}~\cite{Krishnapriyan}.

We use a DNN that has 5 hidden layers each with 50 neurons to approximate the solution.
For the standard PINN, we set the number of residual points $N_r=1000$, the number of initial data $N_i=400$, and the number of periodic boundary data $N_{b}=200$. 
In this case, we set the number of Adam iterative steps to $N_{iter}=2000$. 
For the PT-PINN, one pre-training interval  $[0,0.1]$ is considered, and no additional supervised learning part is added, i.e.~$w_{sp}=0$ in Eq.~\eqref{fortrain}. 
The pre-training and formal training steps of the PT-PINN method use the same neural network structure and training data sets as those of the standard PINN.

\begin{figure}[]
\centering
\includegraphics[width=0.75\textwidth]{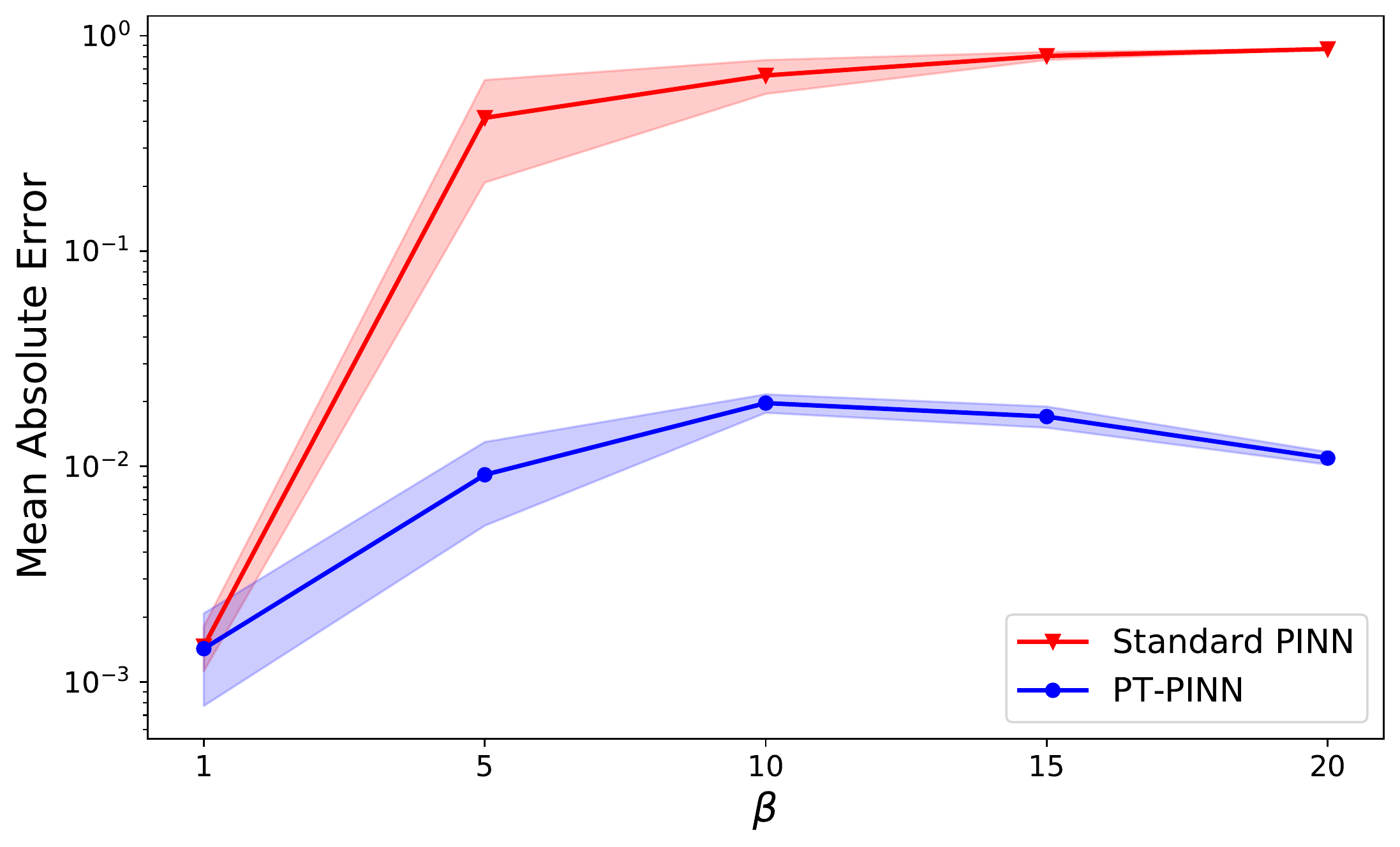}
\caption{Comparison of the absolute errors ($\left\|e\right\|_1$) between the standard PINN and the PT-PINN for Eq.~\eqref{rea} with the different reaction coefficient $\rho=1,5,\cdots,20$. The line and shaded region represent the mean and the standard deviation of 5 independent runs, respectively.}
\label{fig:expN21}
\end{figure}

In this case, we use the the mean absolute errors $\left\|e\right\|_1$ to test the performance of PT-PINN method, which is the same as Refs.~\cite{Krishnapriyan, Haitsiukevich}. We plot the training results of the standard PINN and the PT-PINN method with logarithmic scale in Figure \ref{fig:expN21}. 
For $\rho=1$, both two methods have good performance. For $\rho\geq5$, the standard PINN cannot provide credible prediction any more, while the PT-PINN still gives accurate prediction.

\begin{figure}[]
\centering
\includegraphics[width=0.9\textwidth]{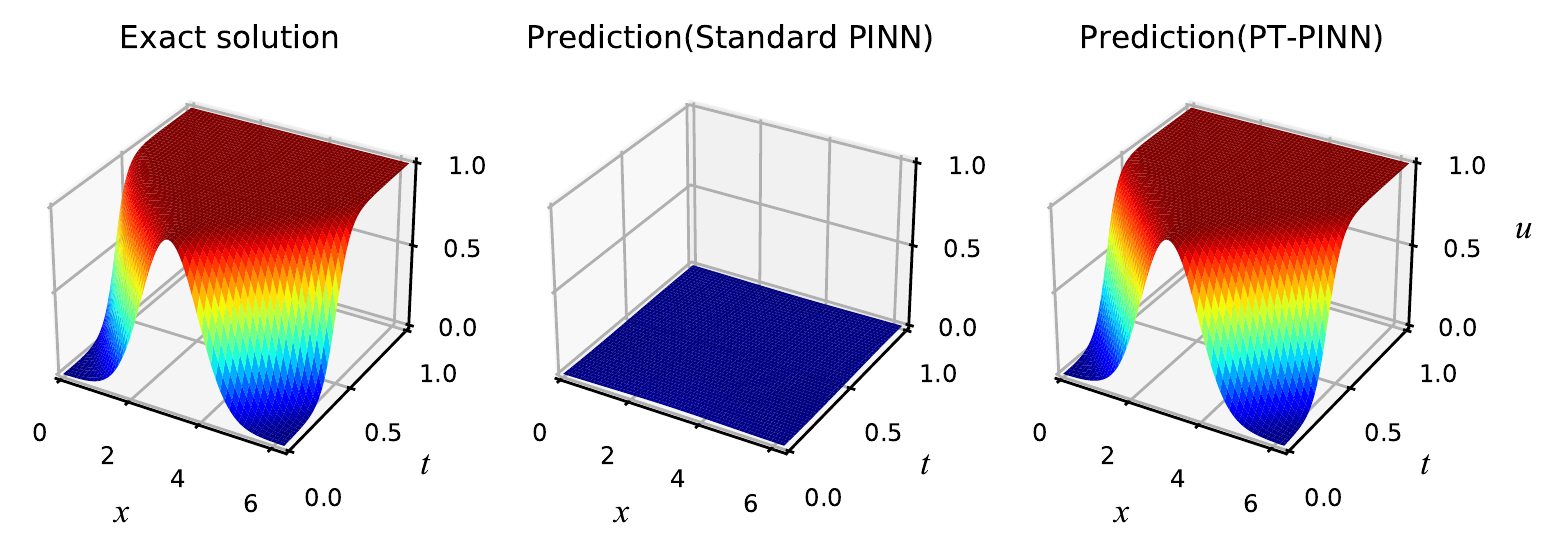}
\caption{The predictions using two PINN methods and the exact solution of Eq.~\eqref{rea} at $\rho=20$.}
\label{fig:expN22}
\end{figure}

\begin{figure}[]
\centering
\includegraphics[width=0.8\textwidth]{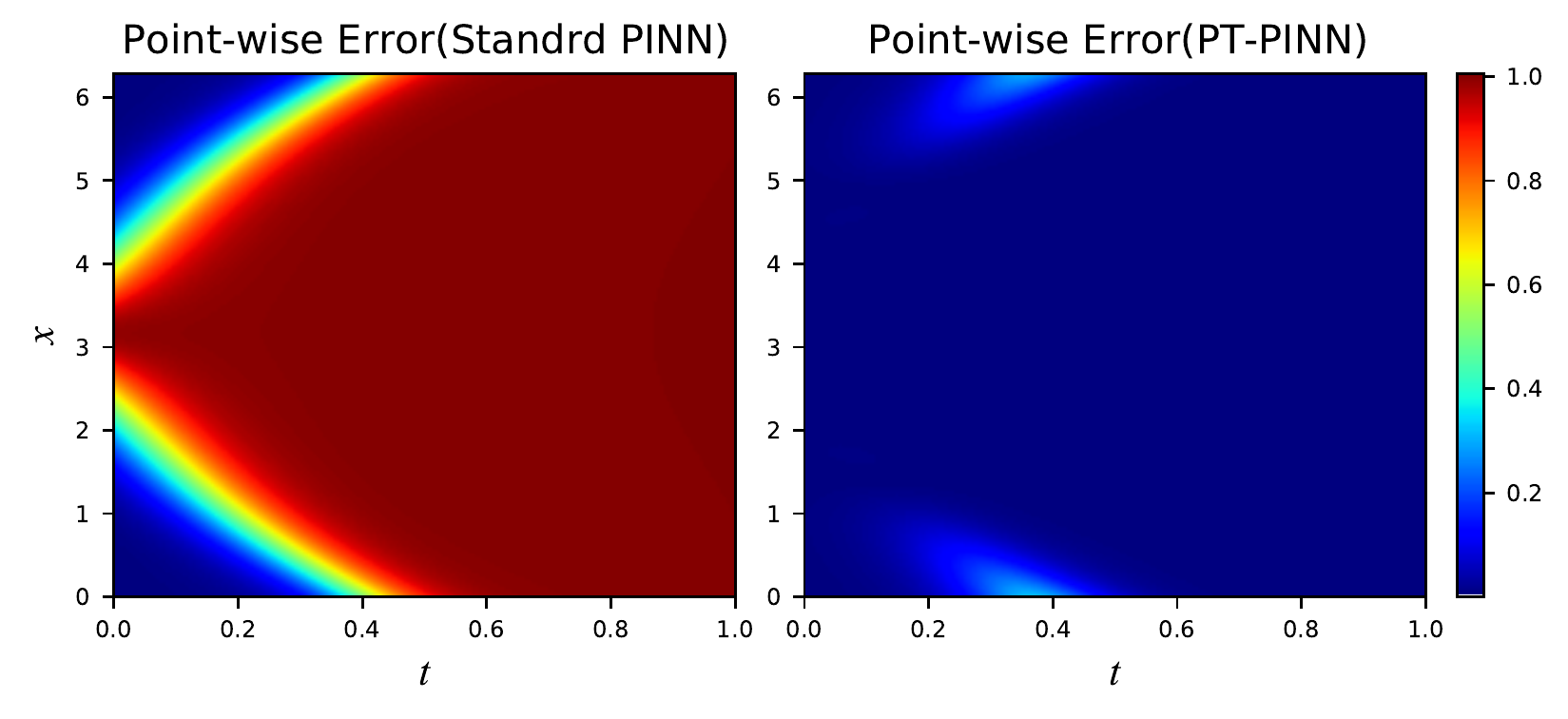}
\caption{Point-wise absolute errors using two PINN methods for Eq.~\eqref{rea} at $\rho=20$.}
\label{fig:expN23}
\end{figure}

We also plot the predictions and the point-wise absolute errors by the two PINN methods for $\rho=20$ in Figure \ref{fig:expN22} and Figure \ref{fig:expN23}, respectively. 
What we observe from Figure~\ref{fig:expN22} is that the standard PINN gives a prediction that is close to zeros in the whole domain, which confirms the results reported in Ref.~\cite{Krishnapriyan}.
The PT-PINN method significantly improves the training performance because of a good initialization of the network parameters provided by the pre-training step.

Table~\ref{tab:exp21} presents the training
errors by using three PINN methods (the standard PINN, the curriculum training method Ref.~\cite{Krishnapriyan} and the 1-interval PT-PINN method. As shown in the table, two methods can provide an effective prediction when $\rho<10$.
Note that, the number of time segments in Refs.~\cite{Krishnapriyan} is 20, which means that it needs to train 20 times to obtain the prediction in whole state space.
The PT-PINN method only trains the networks one more time than the standard PINN method, and does not need to tune the hyper-parameter extensively for each training step.
Moreover, the formal training step of PT-PINN and standard PINN have the same loss function, and the only difference between the two methods are the generations of initial guess. 
We can learn from this example that a bad initial guess will result in training failure.

\begin{table}[]
		\setlength{\abovecaptionskip}{0cm}
		\setlength{\belowcaptionskip}{0.2cm}
	\centering
	\caption{\label{tab:exp21} Means and standard deviations (obtained by 5 independently repeated experiments) of the mean absolute errors ($\|e\|_1$) using the standard PINN method, the ``Time Marching'' method \cite{Krishnapriyan}, and the 1-interval PT-PINN (No extra supervised) for the reaction system Eq.~\eqref{rea} with different reaction coefficient~$\rho$.
	 }
	\begin{tabular}{c|c|c|c}
		\hline
		$\rho$ & Standard PINN & Time Marching method \cite{Krishnapriyan}&  PT-PINN \\ 
		 \hline
		  5 &    $4.2\pm2.1\times10^{-1}$  &  $2.39\times10^{-2}$ &  $9.1\pm3.8\times10^{-3}$ \\
	  	 10 &    $6.5\pm1.2\times10^{-1}$  &  $2.85\times10^{-2}$ &  $2.0\pm0.2\times10^{-2}$ \\
		 15 &    $8.1\pm0.3\times10^{-1}$  &        --            &  $1.7\pm0.2\times10^{-2}$ \\
		 20 &    $8.7\pm0.01\times10^{-1}$ &        --  		  &  $1.1\pm0.1\times10^{-2}$ \\
		 \hline
	\end{tabular}
\end{table}

\subsection{2D heat equation with high-frequency solution}\label{S42}


We consider the 2D linear heat equation as follows:
\begin{equation}\label{Eq42}
\begin{cases}
u_t - (u_{xx}+u_{yy}) = Q(x,y,t),(x,y)\in(0,1)^2,&t\in(0,1],\\
u(x,y,0)=2+\sin(\pi xy),&(x,y)\in(0,1)^2,\\
u(0,y,t)=2+\sin(30\pi t),&t\in(0,1],\\
u(1,y,t)=2+\sin(30\pi t+\pi y)),&t\in(0,1],\\
u(x,0,t)=2+\sin(30\pi t),&t\in(0,1],\\
u(x,1,t)=2+\sin(30\pi t+\pi x),&t\in(0,1].
\end{cases}
\end{equation}
We investigate the exact solution $u(x,y,t)=2+\sin(30 \pi t+\pi xy)$ and derive the source term $Q(x, y, t)$ from the exact solution.

\begin{table}[]
		\setlength{\abovecaptionskip}{0cm}
		\setlength{\belowcaptionskip}{0.2cm}
	\centering
	\caption{\label{tab:exp42} Summary of the training sets of the PINN and PT-PINN methods in Sec.~\ref{S42}} 
\begin{tabular}{c|c|ccc}
\hline
\multicolumn{2}{c|}{\multirow{1}{*}{Method}}&  Time domain & $N_{iter}$ & $~N_r,~N_b,~N_i,~N_{sp}$\\
\hline
\multicolumn{2}{c|}{\multirow{1}{*}{Standard PINN}}&$[0,1.0]$ &    $5000$&  $4000,800,400,-~~$\\
\hline
\multirow{2}{*}{PT-PINN} & 
\multirow{1}{*}{Pre-training}&$[0,0.1]$ &    $2000$&  $2000,400,400,-~~$\\
\cline{2-2} 
& Formal  &$[0,1.0]$ &    $5000$&  $4000,800,400,200$\\
\hline
\end{tabular}
\end{table}

We use the PT-PINN with a single pre-training interval to solve this case. 
The neural network structure for this case is {ResNet}~\cite{HeKM},  which is comprised of five three-layer-deep residual blocks, and the width of each hidden layer is $50$.
The standard PINN is also used to solve this case.
The training set of two methods are listed in Table \ref{tab:exp42}.

\begin{figure}[]
\centering
\includegraphics[width=0.65\textwidth]{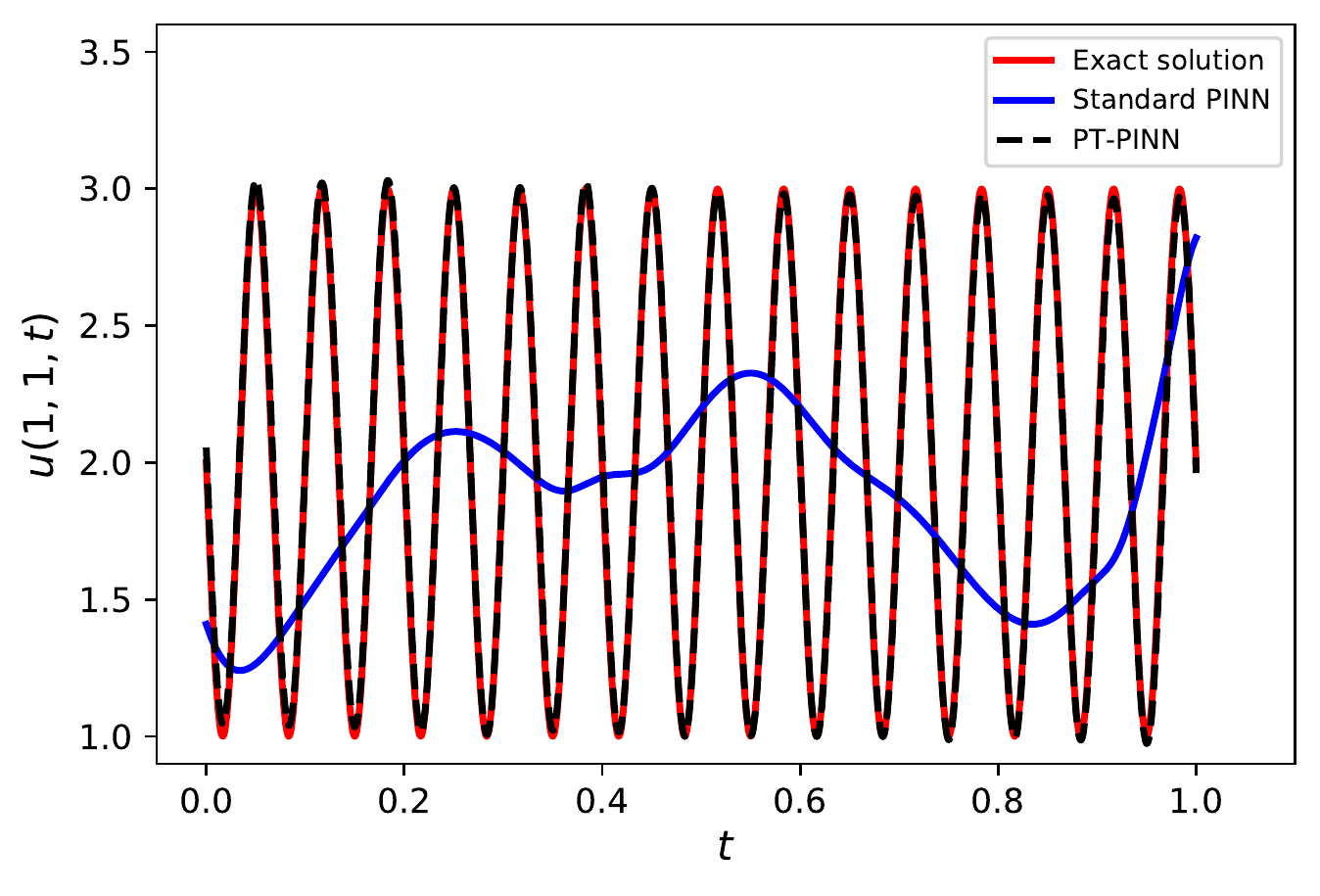}
\caption{The exact solution $u(1,1,t)$ and predictions $u_\theta(1,1,t)$ using the standard PINN and the PT-PINN for Eq.~\eqref{Eq42}
}
\label{fig:exp42}
\end{figure}

Figure \ref{fig:exp42} presents a comparison between the exact solution and the predictions of the two PINN methods at the point $(x,y)=(1,1)$.
It is observed that the prediction of the PT-PINN matches well with the exact solution, whereas the standard PINN cannot correctly solve this equation. 
Table~\ref{tab:exp421} shows the numerical errors of the standard and the PT-PINN methods in the whole spatial-temporal domain.

\begin{table}[]
		\setlength{\abovecaptionskip}{0cm}
		\setlength{\belowcaptionskip}{0.2cm}
	\centering
	\caption{\label{tab:exp421} Test errors of two PINN methods for Eq.~\eqref{Eq42}}
	\begin{tabular}{c|ccc}
		\hline
		Method  &  $\left\|\epsilon\right\|_2$ & $\left\|e\right\|_1$ & $\left\|e\right\|_\infty$ \\ 
		 \hline
	Standard PINN   &  $4.849852\times10^{-1}$ & $8.344807\times10^{-1}$ & $3.990169\times10^0~~$ \\
	PT-PINN&  $3.513617\times10^{-3}$ & $ 5.646592\times10^{-3}$ & $4.673612\times10^{-2}$ \\
		\hline
	\end{tabular}
\end{table}



\subsection{Heat equations with strongly nonlinear diffusion coefficients}\label{S43}

We consider a nonlinear heat equation defined by
\begin{equation}\label{exp43}
\begin{cases}
\displaystyle\frac{\partial u}{\partial t}-\nabla\cdot \left(u^l \nabla u\right)=Q\left(x,t\right), &x\in(0,1),t\in(0,1],\\
u(x,0)=1+2\sin(\pi x),&x\in(0,1),\\
u(0,t)=1+2\sin(2\pi t),&t\in(0,1],\\
u(1,t)=1+2\sin(2\pi t+\pi),&t\in(0,1].
\end{cases}
\end{equation}
In this case, we use an exact solution of $u(x,t)=1+2\sin(\pi (2t+x))$, and derive the source term $Q$ by substituting the exact solution in the equation.
We set the parameter $l=0,1,\cdots,4$ to introduce different degree of non-linearity in the equation, and test the performance of the PT-PINN method.

The neural network structure in this case is the same as that reported in Sect.~\ref{reaction}. 
We set the number of residual points $N_r=1000$, the numbers of supervised points $N_i=N_b=400$ and the number of Adam iterations $N_{iter}=5000$. 
For the PT-PINN, we use one pre-training interval $[0,0.1]$.
No additional supervised learning part is added in the formal training .
The PT-PINN method and the standard PINN method share the same neural network structure and the same number of training data.

\begin{figure}[]
\centering
\includegraphics[width=0.65\textwidth]{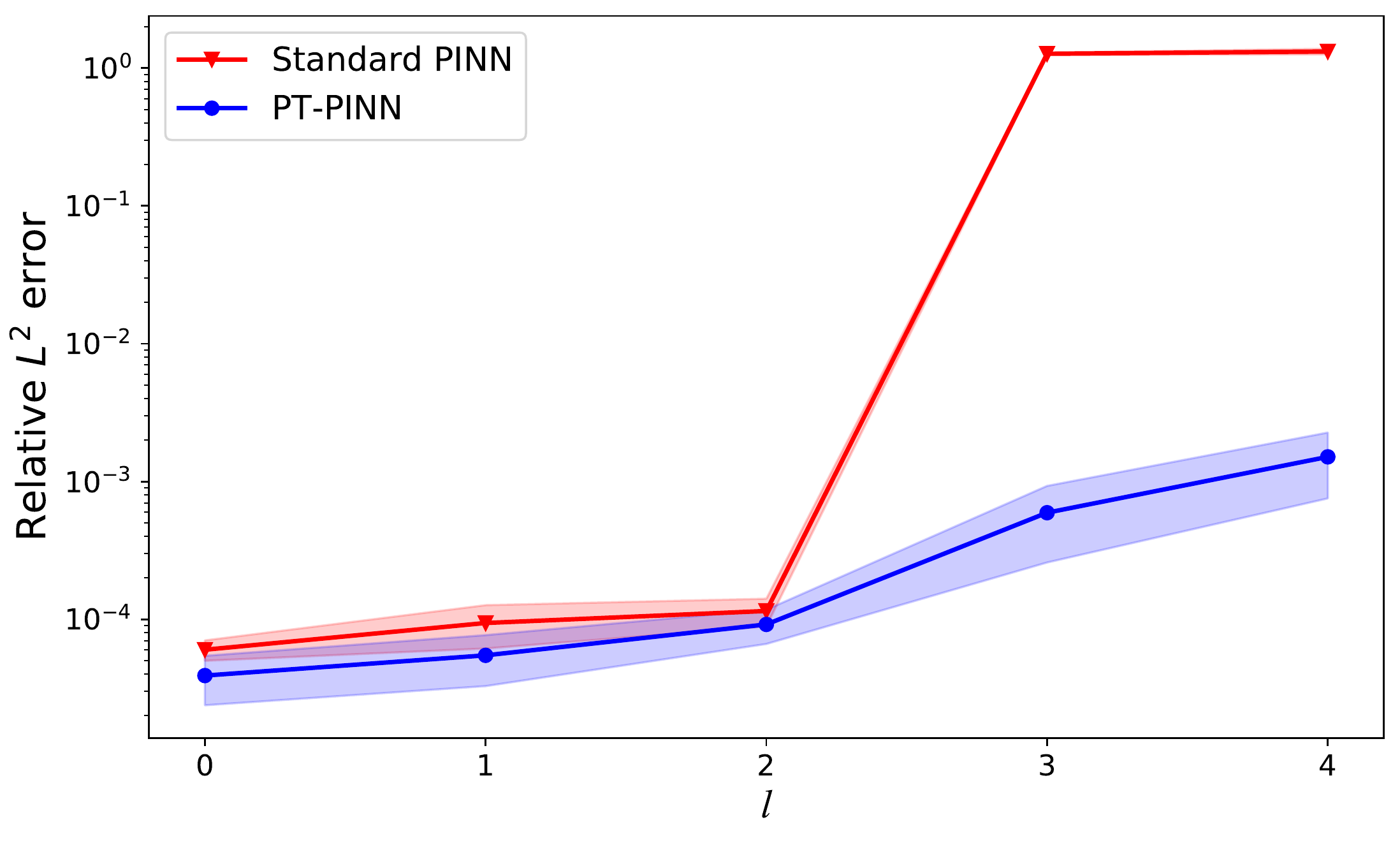}
\caption{Comparison of $L^2$ relative error between the standard PINN and the PT-PINN for Eq.~\eqref{exp43} with different diffusion coefficient $\kappa(u)=u^l, (l=0,1,\cdots,4)$. The line and shaded region represent the mean and the standard deviation of 5 independent runs.}
\label{fig:exp43}
\end{figure}

We plot the relative errors of the standard PINN and the PT-PINN method as a fcuntion of increasing non-linearity (increasing $l$) in Figure \ref{fig:exp43}. 
For relatively weak non-linearity ($l=0,1,2$), both two methods have good performance, and the error in the prediction of the PT-PINN overlaps with the standard PINN within in the standard deviation bar.
For relatively strong non-linearity ($l=3,4$), 
the standard PINN cannot provide credible prediction any more, while the PT-PINN using the pre-training strategy still gives accurate prediction. 
The relative $L^2$ error of both methods are listed in Table~\ref{tab:exp43}, and each error is an average of 5 independently repeated experiments.

\begin{table}[]
		\setlength{\abovecaptionskip}{0cm}
		\setlength{\belowcaptionskip}{0.2cm}
	\centering
	\caption{\label{tab:exp43}Means and standard deviations (obtained by 5 independently repeated experiments) of the relative $L^2$ errors of the two PINN methods for Eq.~\eqref{exp43}} 
	\begin{tabular}{c|c|c}
		\hline
		$l$ & Standard PINN & PT-PINN \\ 
		 \hline
		 0 &    $6.02\pm1.02\times10^{-5}$  &     $3.91\pm1.53\times10^{-5}$ \\
	  	 1 &    $9.41\pm3.27\times10^{-5}$  &     $5.48\pm2.20\times10^{-5}$ \\
		 2 &    $1.15\pm0.26\times10^{-4}$  &     $9.19\pm2.56\times10^{-5}$  \\
		 3 &    $1.27\pm0.03\times10^{0}~~$   &     $5.94\pm3.35\times10^{-4}$ \\
		 4 &    $1.32\pm0.06\times10^{0}~~$    &    $1.51\pm0.76\times10^{-3}$\\
		\hline
	\end{tabular}
\end{table}

\subsection{Allen-Cahn equation with strongly nonlinear source term}\label{S44}

The Allen-Cahn(AC) equation is a widely used PDE for studying the phenomena of phase separation~\cite{Bazant,Allen}. 
In this work, the following spatially 1D time varying AC equation is considered
\begin{equation}\label{AC}
\begin{cases}
u_t-0.0001u_{xx} = 5\left(u^3-u\right),&x\in(-1,1),t\in(0,1],\\
u(x,0)=x^2\cos(\pi x),&x\in(-1,1),\\
u(-1,t)=u(1,t),~u_x(-1,t)=u_x(1,t),&t\in(0,1].
\end{cases}
\end{equation}
The reference solution of Eq. \eqref{AC} is provided by Ref. \cite{RAISSI1}. Many researchers try to solve this AC equation by using PINN method, and the standard PINN failed to provide an efficient prediction~\cite{RAISSI1,Mattey,Wight,Wang4}.
In Ref. \cite{RAISSI1}, the authors propose a time-discretized PINN method by the aid of Runge-Kutta algorithm. 
In Ref.~\cite{Mattey}, it is pointed out that the main reason of the failure of standard PINN in AC equation is that it fails to predict the non-linear source term $5(u^3-u)$ (see Eq.~\eqref{AC}).

We first solve the AC equation by the standard PINN to illustrate difficulties encountered in the training process, and then use 1-interval and 2-interval PT-PINN to overcome these difficulties. 
In this case, we adopt fully-connected DNN having 5 hidden layers each with 50 neurons.
The training data sets used in the example are listed in Table~\ref{tab:AC}. Note that, the 2-interval PT-PINN is designed by performing one more pre-training step on the basis of the 1-interval PT-PINN.

\begin{table}[]
		\setlength{\abovecaptionskip}{0cm}
		\setlength{\belowcaptionskip}{0.2cm}
	\centering
	\caption{\label{tab:AC} The training sets of the PINN, the 1-interval PT-PINN and the 2-interval PT-PINN for the AC equation Eq.~\eqref{AC}.} 
\begin{tabular}{c|c|cc}
\hline
\multicolumn{2}{c|}{\multirow{1}{*}{Method}}&  Time\! domain & $N_r,N_b,N_i,N_{sp}$\\
\hline
\multicolumn{2}{c|}{\multirow{1}{*}{Standard PINN}}&$[0,1.0]$ &  $4000,200,400,-~~~~$\\
\hline
\multirow{2}{*}{\makecell[c]{PT-PINN \\ (1-interval)}} & 
\multirow{1}{*}{1st pre-training}&$[0,0.25]$ &   $2000,200,400,-~~~~$\\
\cline{2-2} 
 &Formal training &$[0,1.0]$ &   $4000,200,400,1000$\\
\hline
\multirow{3}{*}{\makecell[c]{PT-PINN \\ (2-interval)}} & 
\multirow{1}{*}{1st pre-training}&$[0,0.25]$ &   $2000,200,400,-~~~~$\\
\cline{2-2} 
 &2nd pre-training  &$[0,0.5]$ &    $3000,200,400,1000$\\
 \cline{2-2} 
 &Formal training &$[0,1.0]$ &    $4000,200,400,1000$\\
\hline
\end{tabular}
\end{table}

\begin{figure}[]
\centering
\includegraphics[width=0.9\textwidth]{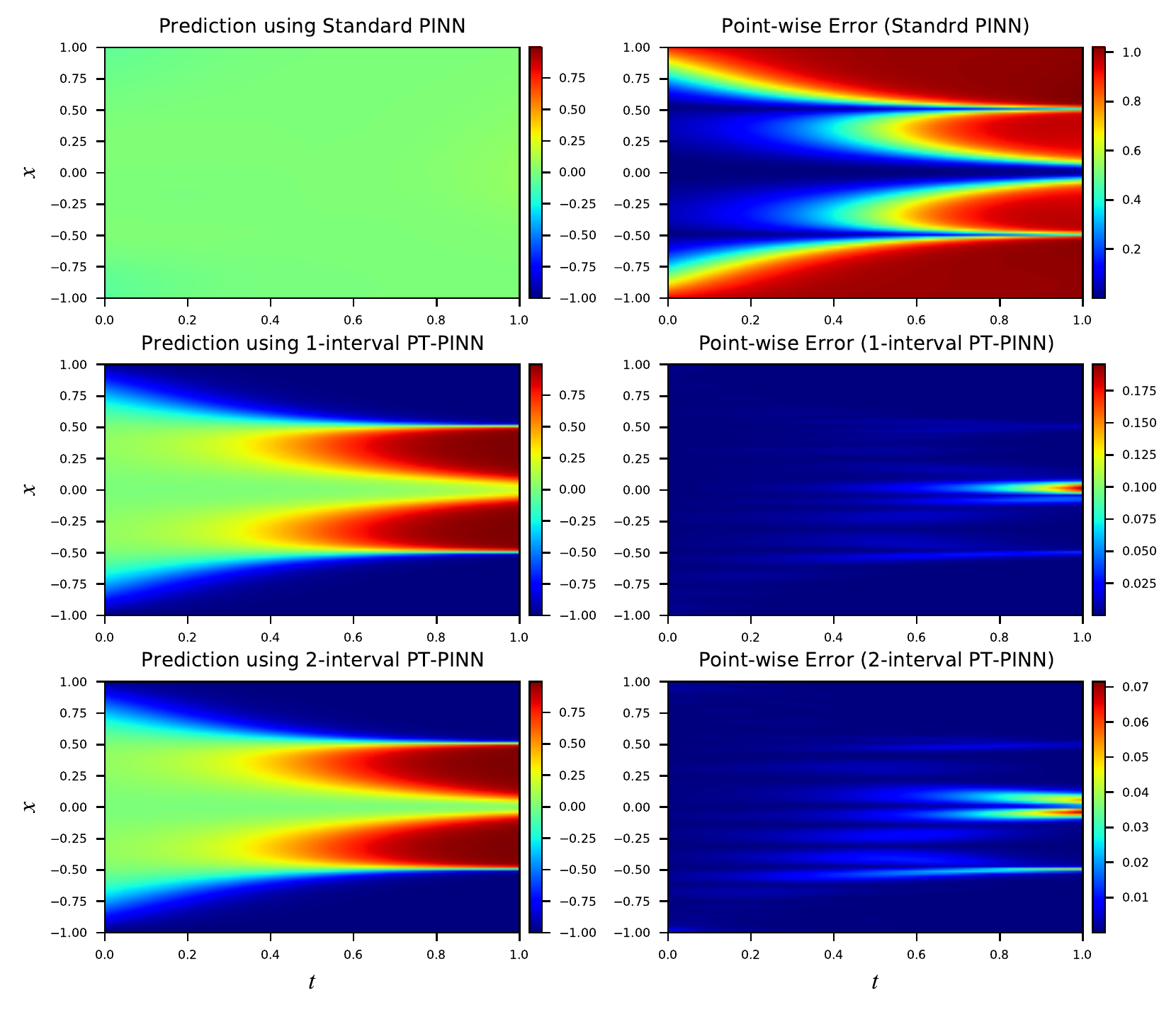}
\caption{The reference solution and the training results (predictions, point-wise absolute errors) obtained by using the standard PINN method and the PT-PINN (1-interval, 2-interval) for  Eq.~\eqref{AC}.}
\label{fig:AC1}
\end{figure}

Figure \ref{fig:AC1} shows the numerical results by using the standard PINN, the 1-interval PT-PINN and the 2-interval PT-PINN.
The standard PINN fails in predicting the solution of the AC equation. 
Both 1-interval and 2-interval PT-PINN provide an accurate approximation to the reference solution of the AC equation.
The 2-interval PT-PINN has a higher accuracy in the regions where the gradient of the solution is large.
\begin{table}[]
		\setlength{\abovecaptionskip}{0cm}
		\setlength{\belowcaptionskip}{0.2cm}
	\centering
	\caption{\label{tab:expac} Means and standard deviations (obtained by 5 independently repeated experiments) of the relative $L^2$ errors ($\left\|\epsilon\right\|_2$) using the standard PINN method, the ``Time Marching'' method \cite{Mattey,Wight}, Self-adaptive PINN \cite{Levi}, and the PT-PINN method for AC equation \eqref{AC}.
	} 
	\begin{tabular}{c|cc}
		\hline
		Method  & $\left\|\varepsilon\right\|_2$ & $N_r$\\ 
		 \hline
		 Standard PINN \cite{RAISSI1}          &    $9.90\times10^{-1}$ & 4000 \\
		 Adaptive resampling \cite{Wight} &      $2.33\times10^{-2}$& - \\
		 Self-adaptive PINN \cite{Levi} & $2.10\times10^{-2}$ & 20000\\
	  	 bc-PINN \cite{Mattey}                &     $1.68\times10^{-2}$&20000 \\
		 PT-PINN(1-interval, this work)     &    ${1.8\pm0.3\times10^{-2}}$&4000\\
		 PT-PINN(2-intervals, this work)    &    $\bm{9.7\pm0.4\times10^{-3}}$&4000\\
		 \hline
	\end{tabular}
\end{table}

Table \ref{tab:expac} provide some results of this example by existing methods in Refs~\cite{RAISSI1,Mattey,Wight,Levi}, and the PT-PINN method with 2-interval is more effective than others. Note that, compared with other methods, the PT-PINN method uses fewer training points and additional training.




\begin{figure}[]
\centering
\includegraphics[width=0.9\textwidth]{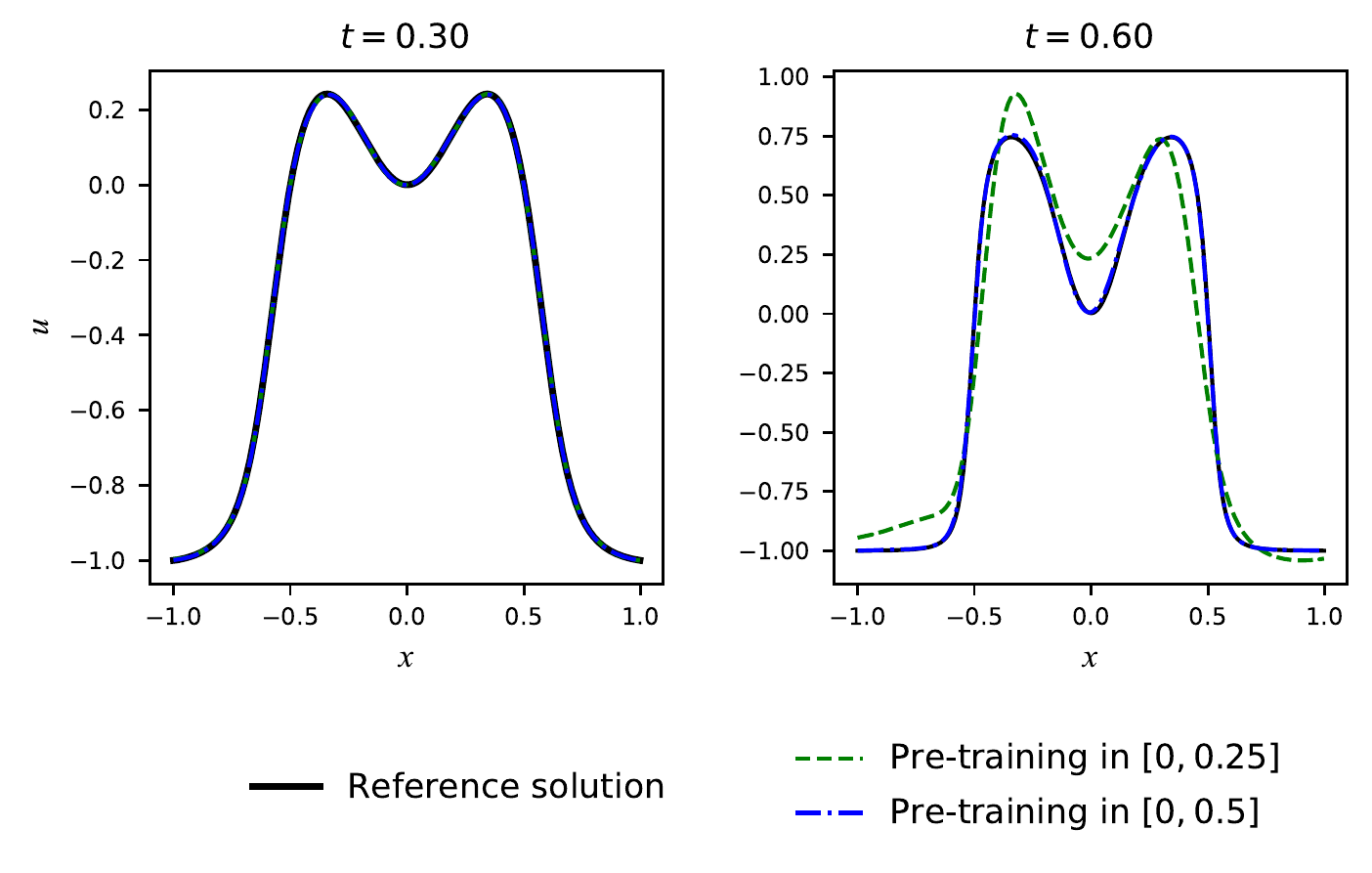}
\caption{The reference solution and the predictions of two pre-training models in 2-interval PT-PINN at different time $(t=0.3,0.6)$ for AC Eq.~\eqref{AC}}
\label{fig:AC3}
\end{figure}

Figure \ref{fig:AC3} shows the comparison of the reference solution and the predictions of the pre-trained models in different time intervals, i.e.~models pre-trained in $[0,0.25]$ and $[0,0.5]$ intervals.
The time domain of the first pre-training is $[0,0.25]$.
Its prediction at $t=0.3$, a time not far away from the pre-training interval, matches the reference solution very well. 
At $t=0.6$, the accuracy of the first pre-training significantly decreases, but it still predicts a solution that roughly has the same shape as the reference.
The time domain of the second pre-training is $[0,0.5]$. Its prediction at $t=0.6$ accurately matches the reference solution. 
The results of the second pre-training can be directly used to help the formal training to obtain a reasonably accurate approximation of the solution in $t\in [0,1]$, as the 2-interval PT-PINN method shows in Figure \ref{fig:AC1}.
It should be pointed out that the 1-interval PT-PINN method can not give an ideal prediction at the positions where the gradient of the solution is large, and the accuracy can be improved obviously by using 2-interval PT-PINN.

\subsection{Convection equation with large convection coefficient}\label{Con}

In this case, we consider using the PT-PINN method to solve a spatially 1D convection equation that has the following form Refs.~\cite{Krishnapriyan,Haitsiukevich}:
\begin{equation}\label{1Dconvection}
\begin{cases}
\displaystyle\frac{\partial u}{\partial t} + \beta \frac{\partial u}{\partial x} = 0,&(x,t)\in(0,2\pi)\times(0,1],\\
u(x,0) = \sin(x),&x\in(0,2\pi),\\
u(0,t)=u(2\pi,t),&t\in(0,1].
\end{cases}
\end{equation}
The exact solution is $u(x,t)=\sin(x-\beta t)$. Ref.~\cite{Krishnapriyan} shows that the standard PINN method fails to learn the solution in the cases of relatively large convection coefficients $\beta\geq 10$, and the author improves the
performance of PINN with relatively large $\beta$ value $30\leq \beta \leq 40$ based on the ``curriculum training'' method.

\begin{table}[]
		\setlength{\abovecaptionskip}{0cm}
		\setlength{\belowcaptionskip}{0.2cm}
	\centering
	\caption{\label{tab:exp61} Means and standard deviations (obtained by 5 independently repeated experiments) of the mean absolute errors ($\|e\|_1$) using the standard PINN method, the curriculum training method \cite{Krishnapriyan}, and the 2-interval PT-PINN  for Eq.~\eqref{1Dconvection} with
	different convection coefficient~$\beta$.
	$\ast$: At $\beta=60$, 10 independent experiments were conducted. One of them failed (see the main text for discussion). The average and standard deviation are taken over the 9 successful cases.
	} 
	\begin{tabular}{c|c|c|c}
		\hline
		$\beta$ & Standard PINN & Curriculum training \cite{Krishnapriyan} & PT-PINN \\ 
		 \hline
		 20 &    $1.9\pm1.0\times10^{-3}$  &  $5.42\times10^{-3}$&   $6.3\pm0.3\times10^{-4}$ \\
	  	 30 &    $2.8\pm1.0\times10^{-3}$  &   $1.10\times10^{-2}$&    $1.8\pm0.5\times10^{-3}$ \\
		 40 &    $2.0\pm1.6\times10^{-1}$  & $2.69\times10^{-2}$&       $2.3\pm0.4\times10^{-3}$  \\
		 50 &    $3.7\pm0.3\times10^{-1}$ & --&    $4.8\pm1.3\times10^{-3}$ \\
		 60* &    $4.4\pm0.2\times10^{-1}$ & --&    $4.9\pm2.6\times10^{-3}$ \\
		 \hline
	\end{tabular}
\end{table}

\begin{figure}[]
\centering
\includegraphics[width=1\textwidth]{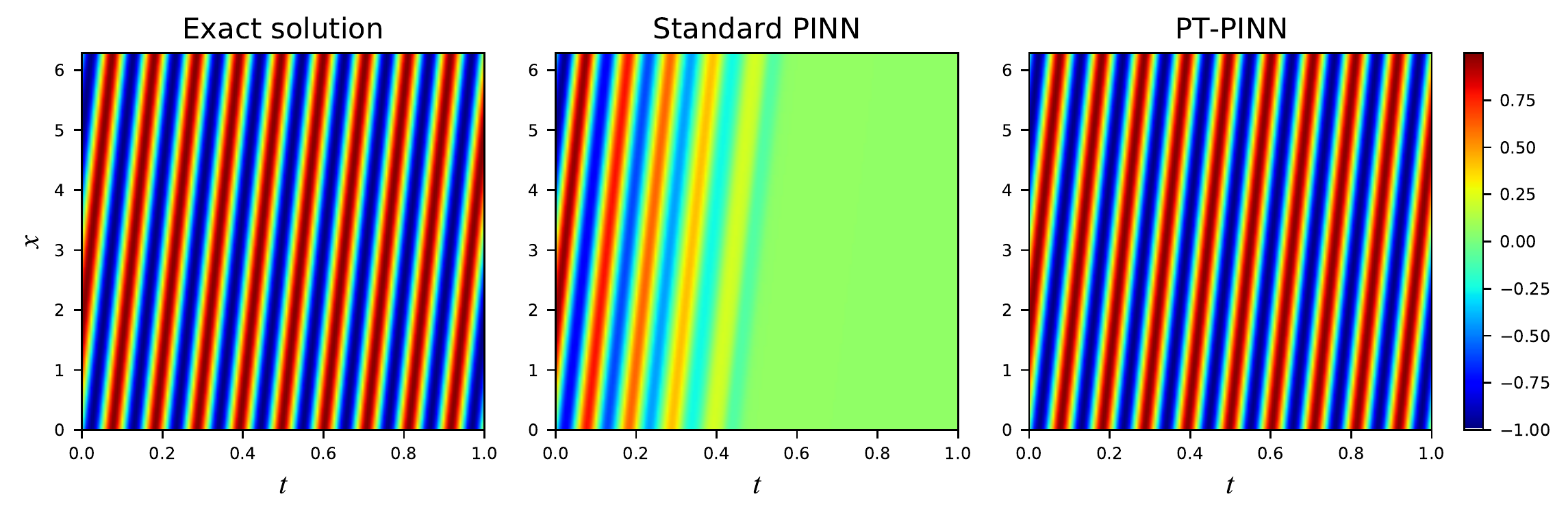}
\caption{Exact solution and predictions using two PINN methods for Eq.~\eqref{1Dconvection} at $\beta=60$.}
\label{Figcon1}
\end{figure}

In this work, the neural network structure is ResNet which is  same as Sect.~\ref{S42}. 
The numbers of training data are set to $N_i=400,N_b=400,N_r=2000$.
The same size of training data set is used for the standard PINN and the pre-training and formal training steps of the PT-PINN. 
The number of Adam iterations is $N_{iter}=2000$.
The PT-PINN adopts two pre-training intervals: $[0,0.2]$ and $[0,0.6]$, and the number of extra supervised learning points is set to $N_{sp}=1000$. 

For comparing the performances of different methods, we use the mean absolute error $\left\|e\right\|_1$ to evaluate the training results.
Table \ref{tab:exp61} presents the training errors by using three PINN methods (the standard PINN, the curriculum training method~\cite{Krishnapriyan} and the PT-PINN). 
The training of the convection equation becomes more difficult as the convection coefficient increases. 
As shown in the table, the standard PINN  can provide an accurate prediction with relatively small convection coefficients ($\beta=20,30$).
In the case of $\beta=40$, the standard PINN fails to work, but the curriculum training provides an effective prediction.
The errors of the curriculum training  at larger convection coefficients are not reported by Ref.~\cite{Krishnapriyan}.
The PT-PINN succeeds in solving the convection equation
in the whole range of convection coefficient $20\leq \beta \leq 60$, and its errors are the lowest among the three compared methods. The PT-PINN achieves satisfactory accuracy even in the most difficult case of $\beta = 60$, and is compared with the standard PINN in Figure~\ref{Figcon1}. 
In Ref.~\cite{Krishnapriyan}, the author shows that large convection coefficient results in a complex and non-symmetric loss landscape, and this is the reason why even larger convection coefficients are not considered by this work. Note that the PT-PINN occasionally fail (1 failures out of 10 independent trainings) at $\beta = 60$.
The failed case can be easily picked out because the converged loss is significantly larger than the successful cases, thus we only report the average error of the converged cases at $\beta = 60$ in Table~\ref{tab:exp61}.




\begin{figure}[]
\centering
\includegraphics[width=1\textwidth]{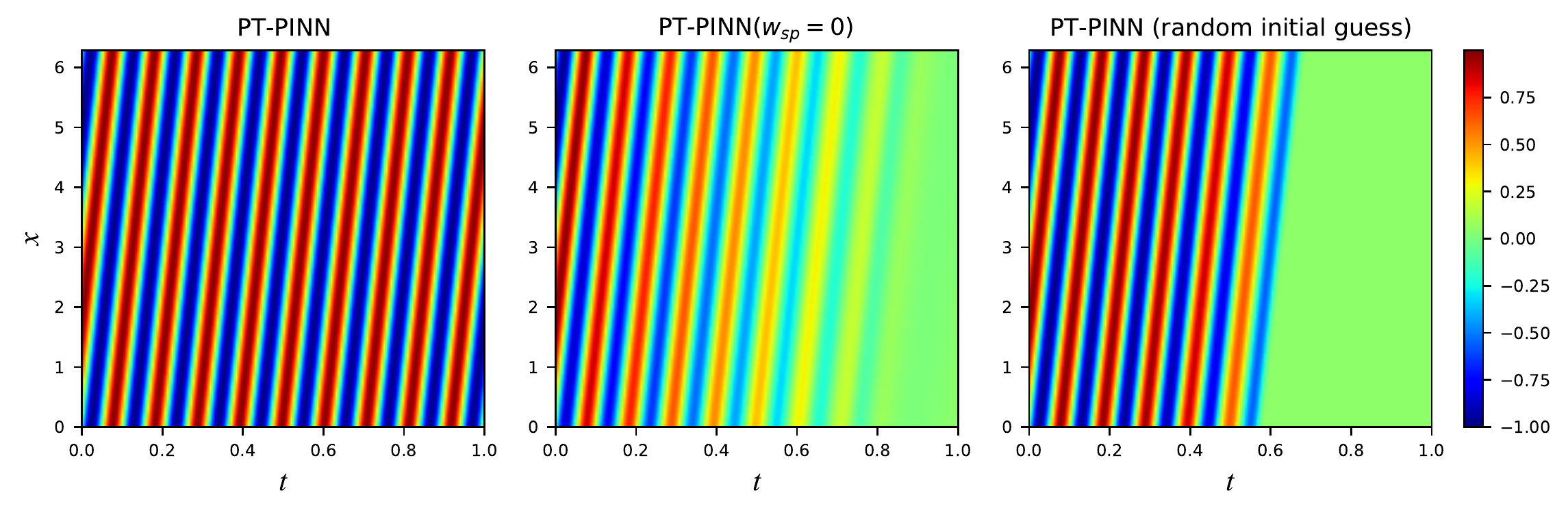}
\caption{The training results using the PT-PINN method and two variants of the PT-PINN ($w_{sp}=0$ and random initial guess) in the formal training step for the convection equation Eq.~\eqref{1Dconvection} at $\beta=60$.
}
\label{Figcon12}
\end{figure}

The PT-PINN improves the performance of training by using two kinds of information from the pre-training step: 
the extra supervised learning data and the pre-trained neural network parameters for the network initialization in the formal training step. 
We present an ablation study in Figure~\ref{Figcon12} to demonstrate the effect of these two kinds of information.
The PT-PINN is presented in the left plot.
In the middle plot, the solution using no supervised training data ($w_{sp} = 0$) in the formal training step is presented. 
It has a periodic structure, but the frequency and the magnitude at large $t$ deviate from those of the exact solution.
The right plot gives the solution trained with a random network parameter initialization in the formal training step.
It fails to capture the periodic structure at time $t > 0.5$.
It is thus clear from the Figure~\ref{Figcon12} that both the addition supervised learning data and the pre-trained initialization of the network parameters are crucial for solving the convection equation with a high-frequency solution.

\section{Conclusion}
In this work we propose the pre-training physical informed neural network (PT-PINN) method for solving the evolution PDEs.
The multiple pre-training strategy transforms a difficult problem into several relatively simple problems by shortening the time domain length, so as to provide a good initialization of the neural network parameters and an extra supervised learning data set for the formal training step on the entire time domain.
The PT-PINN can dramatically improve the training performance by pre-training strategy. 
We demonstrate that both the initialization from pre-training and the supervised learning data are crucial for improving the training performance.
By numerical examples, we show that the standard PINN may fail in solving evolution PDEs with strong non-linearity or having a high-frequency solution, while the proposed PT-PINN method succeeds in almost all the testing cases and has a substantially higher accuracy than the standard PINN method.
When compared with the state-of-the-art time marching and the curriculum regularization methods, the PT-PINN achieves higher accuracy with fewer residual points and fewer training costs.



It is noted that the success of the PT-PINN relies on the assumption that the solution of the PDE is smooth, thus it may be extrapolated out of the training time-domain to a certain extend (see detailed discussion in Sect.~\ref{S44}).
It is likely that the PT-PINN method would fail in the cases that the solutions present discontinuity. In this work, the length of pre-training intervals have to be specified in advance, and according to our experience, the different lengths of intervals will affect the training performance to some extent.
How to design an effective and adaptive pre-training strategy for the evolution PDEs with non-smooth solutions is an open problem that will be investigated in the future. 

\section*{Code availability}

The code of this work is publicly available online via  \url{https://doi.org/10.5281/zenodo.6642181}.


\section*{Acknowledgements}
The work of J.G.~was supported by the Project of Natural Science Foundation of Shandong Province No.ZR2021MA092.
The work of Y.Y.~was supported by the foundation of CAEP (CX20210044). 
The work of H.W.~was supported by the National Science Foundation of China under Grant No.11871110 and 12122103. 
The work of T.G.~was supported by the laboratory of computational physics.

\bibliography{mybibfile}

\end{document}